\documentclass[12pt]{amsart}
\usepackage{amssymb}

\begin{document}
\baselineskip 17.2pt

\newtheorem{theorem}{Theorem}[section]
\newtheorem{prop}[theorem]{Proposition}
\newtheorem{lemma}[theorem]{Lemma}
\newtheorem{definition}[theorem]{Definition}
\newtheorem{corollary}[theorem]{Corollary}
\newtheorem{example}[theorem]{Example}
\newtheorem{remark}[theorem]{Remark}
\newcommand{\ra}{\rightarrow}
\renewcommand{\theequation}
{\thesection.\arabic{equation}}
\newcommand{\ccc}{{\cal C}}

\def \LVR {${\rm  L}^{p, \lambda}({\mathbb R}^n)$}
\def \LNVR {${\rm  L}_L^{p, \lambda}({\mathbb R}^n)$}
\def \LV {${\rm  L}^{p, \lambda}({\mathbb R}^n)$}
\def \LNV {${\rm  L}_L^{p, \lambda}({\mathbb R}^n)$}
\def \LM {${\mathcal  M}^{p, \lambda}(\Omega)$}
\def \LNVO {${\mathcal  L}_L^{p, \lambda}(\Omega)$}
\def \LNO {${\mathcal  L}_L^{p', \lambda}(\Omega)$}
\def\RR{\mathbb R}

\baselineskip 18.00 pt
\title[Old and new  Morrey spaces via heat kernel bounds]
{Old and new Morrey spaces via heat kernel bounds}
\thanks{XTD is supported by a grant from Australia Research
Council; JX is supported by NSERC of Canada; LXY is partially supported
by NSF of China (Grant No. 10371134).
 }

\date{}
\author[X.T. Duong, J. Xiao and L.X. Yan]{ Xuan Thinh Duong,\ \
Jie Xiao\ \  and\ \
Lixin Yan}


\begin{abstract}
 Given $p\in [1,\infty)$ and $\lambda\in (0,n)$, we study Morrey space
 ${\rm L}^{p,\lambda}({\mathbb R}^n)$ of all locally integrable complex-valued functions
 $f$ on ${\mathbb R}^n$ such that for every open Euclidean ball $B\subset{\mathbb R}^n$
 with radius $r_B$ there are numbers $C=C(f)$ (depending on $f$) and $c=c(f,B)$ (relying
 upon $f$ and $B$) satisfying
$$
r_B^{-\lambda}\int_{B}|f(x)-c|^pdx\le C
$$
and derive old and new, two essentially different cases arising from either choosing
$c=f_B=|B|^{-1}\int_{B}f(y)dy$ or replacing $c$ by $P_{t_B}(x)=\int_{t_B}p_{t_B}(x,y)f(y)dy$
 -- where $t_B$ is scaled to $r_B$ and $p_t(\cdot,\cdot)$ is the kernel of the infinitesimal
 generator $L$ of an analytic semigroup $\{e^{-tL}\}_{t\ge0}$ on ${\rm L}^2({\mathbb R}^n)$.
 Consequently, we are led to simultaneously characterize the old and new Morrey spaces, but
 also to show that for a suitable operator $L$, the new Morrey space is equivalent to the old one.
\end{abstract}

\maketitle

\medskip


\section{Introduction}
\setcounter{equation}{0}

As well-known, a priori estimates mixing ${\rm L}^p$ and ${\rm Lip}_\lambda$ are frequently
used in the study of partial differential equations -- naturally, the so-called Morrey spaces
are brought into play (cf. \cite{T}). A locally integrable complex-valued function $f$ on
${\mathbb R}^n$ is said to be in the Morrey space  \LV,  $1\leq p<\infty$ and $\lambda\in (0, n+p)$,
 if for every Euclidean open ball $B\subset{\mathbb R}^n$ with radius $r_B$ there are numbers $C=C(f)$
 (depending on $f$) and $c=c(f,B)$ (relying upon $f$ and $B$) satisfying
$$
r_B^{-\lambda}\int_{B}|f(x)-c|^pdx\le C.
$$
The space of \LV-functions was introduced by Morrey \cite{Mo}. Since then, the space has been
studied extensively -- see for example \cite{C, JN, JTW, P, S, Sta, Z}.

We would like to note that as in \cite{P},
for $1\leq p<\infty$ and $\lambda=n$, the spaces ${\rm L}^{p, n}({\mathbb R}^n)$ are  variants
 of  the classical BMO (bounded mean oscillation) function space of John-Nirenberg. For
 $1\leq p<\infty$ and $\lambda\in (n, n+p)$, the spaces  ${\rm L}^{p, \lambda}({\mathbb R}^n)$
 are variants of the homogeneous Lipschitz spaces ${\rm Lip}_{(\lambda-n)/p}({\mathbb R}^n)$.

Clearly, the remaining cases: $1\leq p<\infty$ and $\lambda \in (0, n)$ are of independent
interest, and hence motivate our investigation. The purpose of this paper is twofold. First,
we explore some new characterizations of \LV\ through the fact that \LV\  consists of all
 locally integrable complex-valued functions $f$ on ${\mathbb R}^n$ satisfying
\begin{equation}
\|f\|_{{\rm L}^{p, \lambda}}=\sup_{B\subset {\mathbb R}^n}
\bigg[r_B^{-\lambda}\int_B|f(x)-f_B|^pdx\bigg]^{1/p} <\infty,
\label{e1.1}
\end{equation}
where the supremum is taken over all Euclidean open balls $B=B(x_0, r_B)$ with center
$x_0$ and radius $r_B$, and $f_B$ stands
for the mean value of $f$ over $B$, i.e.,
\begin{eqnarray*}
f_B= |B|^{-1}\int_Bf(x)dx.
\end{eqnarray*}

The second aim is to use those new characterizations as motives of a continuous study
of \cite{AdX, DY1, DDY, DY3} and so to introduce new Morrey spaces \LNV\ associated with
operators. Roughly speaking,
if $L$ is the
infinitesimal  generator of an analytic semigroup $\{e^{-tL}\}_{t\geq 0}$ on
 ${\rm L}^2(\mathbb{R}^n)$ with
kernel $p_t(x,y)$ which decays fast enough, then we can view $P_tf=
e^{-tL}f$ as an average version of $f$ at the scale $t$ and use the quantity
\begin{eqnarray*}
P_{t_B}f(x)=\int_{{\mathbb R}^n}p_{t_B}(x,y)f(y)dy
\end{eqnarray*}
to replace the mean value $f_B$ in the equivalent semi-norm (\ref{e1.1}) of the original
 Morrey space, where $t_B$ is scaled to the radius of the ball $B.$
Hence we say that a function $f$ (with appropriate bound on its size $|f|$) belongs to the space
\LNV\ (where $1\leq p<\infty$ and $\lambda \in (0, n)$), provided
\begin{eqnarray}
 \|f\|_{{\rm L}^{p, \lambda}_L}=\sup_{B\subset {\mathbb R}^n}\bigg[ r_B^{-\lambda}\int_B
 |f(x)-P_{t_B} f(x)|^p dx\bigg]^{1/p}
 <\infty
\label{e1.2}
\end{eqnarray}
where $t_B=r_B^m$ for a fixed constant $m>0$ -- see the forthcoming Sections 2.2 and 3.1.

We pursue a better understanding of (\ref{e1.1}) and (\ref{e1.2}) through the following aspects:

In Section 2, we collect most useful materials on the bounded holomorphic functional calculus.

In Section 3, we study some  characterizations of \LV\ and \LNV\ and give a criterion
for \LV $\subseteq$ \LNV. The later fact illustrates that \LV\ exists as the minimal Morrey
space, and consequently induces a concept of the maximal Morrey space.

In Section 4, we establish an identity formula associated with the operator $L$. This formula
is a key to handle the quadratic features of the old and new Morrey spaces.

As an immediate continuation of Section 4, Section 5 is devoted to Littlewood-Paley type
 characterizations of \LV\ and ${\rm L}^{p, \lambda}_L({\mathbb R}^n)$ via the predual of
  \LV\ (cf. \cite{Z}) and a number of important estimates for functions in \LV\ and
  ${\rm L}^{p, \lambda}_L({\mathbb R}^n)$. Moreover, we show that for a suitable semigroup
  $\{e^{-tL}\}_{t>0}$, \LNV\ equals \LV\ with equivalent seminorms -- in particular, if $L$
  is either $\triangle$ or $\sqrt{\triangle}$ on ${\mathbb R}^n$, then \LV\ coincides with
  ${\rm L}^{p, \lambda}_{\sqrt{\triangle}}({\mathbb R}^n)$ and
  ${\rm L}^{p, \lambda}_{{\triangle}}({\mathbb R}^n)$.

Throughout, the letters $c, c_1, c_2,...$ will  denote (possibly different)
constants that are independent of the essential variables.

\bigskip

\section{Preliminaries}
\setcounter{equation}{0}

\subsection{ Holomorphic functional calculi of operators.}
 We start with a review of some definitions of
holomorphic functional calculi introduced by McIntosh
\cite{Mc}.
Let $0\leq \omega<\nu<\pi$. We define the closed sector in the complex plane
${\mathbb C}$
$$
S_{\omega}=\{z\in {\mathbb C}: |{\rm arg}z|\leq\omega\}\cup\{0\}
$$
and denote the interior of $S_{\omega}$ by $S_{\omega}^0$.

We employ the following subspaces of the space $H(S_{\nu}^0)$ of all holomorphic
functions on $S_{\nu}^0$:
$$
H_{\infty}(S^0_{\nu}) =\{b\in H(S_{\nu}^0):\ ||b||_{{\infty}}<\infty\},
$$
where
$$
||b||_{\infty}={\rm sup}\{|b(z)|: z\in S^0_{\nu}\}
$$
{and}
$$
\Psi(S^0_{\nu})=\{\psi\in H(S^0_{\nu}): \exists\  \!  s>0,  \
 |\psi(z)|\leq c
|z|^s({1+|z|^{2s}})^{-1}\}.
$$

Given $0\leq\omega<\pi$, a closed operator $L$ in ${\rm L}^2({\mathbb R}^n)$ is
said to be of type $\omega$ if $\sigma(L)\subset S_{\omega}$, and for
each $\nu>\omega,$ there
exists a constant $c_{\nu}$ such that
$$
 \|(L-\lambda {\mathcal I})^{-1}\|_{2,2}=
 \|(L-\lambda {\mathcal I})^{-1}\|_{{\rm L}^2\to{\rm L}^2}\leq c_{\nu}|\lambda|^{-1},
\ \ \ \ \ \lambda\not\in S_{\nu}.
$$
If $L$ is of type $\omega$ and $\psi\in
\Psi(S^0_{\nu})$, we define
$\psi(L)\in {\mathcal L}({\rm L}^2, {\rm L}^2)$ by
\begin{equation}
\psi(L)=\frac{1}{2\pi i}\int_{\Gamma}(L-\lambda
{\mathcal I})^{-1}\psi(\lambda)d\lambda,
\label{e2.1}
\end{equation}
where $\Gamma$ is the contour $\{\xi=re^{\pm i\theta}: r\geq 0\}$ parametrised clockwise
around $S_{\omega}$, and $\omega<\theta<\nu$.   Clearly, this integral is absolutely
 convergent in ${\mathcal  L}({\rm L}^2, {\rm L}^2)$ (which is the class of all bounded
 linear operators on ${\rm L}^2$),
and it is straightforward to show,
 using Cauchy's theorem, that the definition is independent of the choice of
 $\theta\in (\omega, \nu)$.  If, in addition, $L$ is one-one and has dense range
 and if
$b\in H_{\infty}(S^0_{\nu})$, then $b(L)$ can be defined by
\begin{eqnarray*}
b(L)=[\psi(L)]^{-1}(b\psi)(L)\quad\hbox{where}\quad\psi(z) =z(1+z)^{-2}.
\end{eqnarray*}
It can be shown that $b(L)$ is
a well-defined linear operator in
${\rm L}^2({\mathbb R}^n)$.

We say that $L$ has a bounded
$H_{\infty}$ calculus in ${\rm L}^2({\mathbb R}^n) $ provided there exists $c_{\nu, 2}>0$
such that  $b(L)\in {\mathcal  L}({\rm L}^2, {\rm L}^2)$ and
\begin{eqnarray*}
\|b(L)\|_{2,2}=\|b(L)\|_{{\rm L}^2\to {\rm L}^2}\leq c_{\nu, 2}
||b||_{{\infty}}\quad\forall b\in H_{\infty}(S^0_{\nu}).
\end{eqnarray*}
For the conditions and properties of operators which have
holomorphic functional calculi, see \cite{Mc} and \cite{AlDMc} which
also contain a proof of the following   Convergence Lemma.

\medskip

\begin{lemma}\label{l2.1} {Let $X$ be a complex Banach
space. Given $0\le\omega<\nu\le\pi$, let $L$ be an operator of type
$\omega$ which is one-to-one with dense domain and range. Suppose
$\{f_\alpha\}$ is a uniformly bounded net in $H_\infty(S_\nu^0)$,
which converges to $f\in H_\infty(S_\nu^0)$ uniformly on compact
subsets of $S_\nu^0$, such that $\{f_\alpha(L)\}$ is a uniformly
bounded net in the space ${\mathcal L}(X,X)$ of continuous linear
operators on $X$. Then $f(L)\in{\mathcal L}(X,X)$, $f_\alpha(L)u\to
f(L)u$ for all $u\in X$ and
$$
\|f(L)\|=\|f(L)\|_{X\to X}\le\sup_{\alpha}\|f_\alpha(L)\|=\sup_{\alpha}\|f_\alpha(L)\|_{X\to X}.
$$
}
\end{lemma}

 \subsection{Two more assumptions.}\
Let  $L$ be a linear operator of type
$\omega$ on ${\rm L}^2({\mathbb R}^n)$
with $\omega<\pi/2$, hence $L$ generates a holomorphic semigroup
$e^{-zL}$, $0\leq |{\rm Arg}(z)|<\pi/2-\omega$.  Assume the following two
conditions.

\medskip

\noindent
{\bf Assumption (a)}:\ The holomorphic semigroup
$$
\{e^{-zL}\}_{0\leq|{\rm Arg}(z)|<\pi/2-\omega}
$$
is represented by kernel  $p_z(x,y)$
which satisfies an  upper bound
$$
|p_z(x,y)|\leq c_{\theta}h_{|z|}(x,y)\quad\forall x,y\in {\mathbb R}^n
$$
and
$$
|{\rm Arg}(z)|<\pi/2-\theta\quad\hbox{for}\quad \theta>\omega,
$$
where $h_t(\cdot,\cdot)$ is determined by
\begin{equation}
h_t(x,y)=  t^{-{n/m}}g\bigg({\frac{|x-y|}{t^{1/m}}}\bigg),
\label{e2.2}
\end{equation}
in which $m$ is a positive constant and $g$ is a positive, bounded,
decreasing function satisfying
\begin{equation}
\lim\limits_{r\rightarrow\infty}r^{n+\epsilon}g(r)=0\quad\hbox{for\ some}\ \ \epsilon>0.
\label{e2.3}
\end{equation}
\bigskip

\noindent
{\bf Assumption (b)}:\  The operator $L$ has a bounded $H_{\infty}$-calculus in  ${\rm L}^2({\mathbb R}^n)$.

\medskip

Now, we give some  consequences of the assumptions (a)
and (b) which will be used later.

\medskip

First, if
$\{e^{-tL}\}_{t>0}$ is a bounded analytic semigroup on
${\rm L}^2({\mathbb R}^n)$ whose kernel  $p_t(x,y)$ satisfies the
estimates (\ref{e2.2}) and (\ref{e2.3}), then for any $k\in{\mathbb N}$,
the time
derivatives of
$p_t$ satisfy
\begin{eqnarray}
\Big|t^k{\frac{\partial^k p_t(x,y)}{\partial t^k}}\Big|\leq {\frac{c}{t^{n/m}}}
 g\Big(\frac{|x-y|}{t^{1/m}}\Big)\quad\hbox{for\ all}\ t>0\ \ \hbox{and\ almost\ all}\  x,y\in{\mathbb R}^n.
\label{e2.4}
\end{eqnarray}
For each   $k\in{\mathbb N}$, the function $g$ might depend on $k$ but it always
  satisfies (\ref{e2.3}). See
 Theorem 6.17  of \cite{O}.

\medskip

Secondly, $L$ has a bounded $H_{\infty}$-calculus in  ${\rm L}^2({\mathbb R}^n)$  if and only if
for any non-zero function $\psi\in \Psi(S^0_{\nu})$, $L$  satisfies
the square function estimate    and its reverse
\begin{eqnarray}
c_1\|f\|_{{\rm L}^2}\leq \Big( \int_0^{\infty}\|\psi_t(L)f\|_{{\rm L}^2}^2
{\frac{dt}{t}}\Big)^{1/2}\leq c_2\|f\|_{{\rm L}^2}
\label{e2.5}
\end{eqnarray}
 for some $0<c_1\leq c_2<\infty$,
where $\psi_t(\xi)=\psi(t\xi)$.
Note that different choices of $\nu>\omega$ and $\psi\in \Psi(S^0_{\nu})$
lead to equivalent quadratic norms of $f.$

As noted in \cite{Mc}, positive self-adjoint operators satisfy the
quadratic estimate (\ref{e2.5}). So do normal operators with spectra
in a sector, and maximal accretive operators. For definitions of
these classes of operators, we refer the reader to \cite{Y}.

The following result, existing as a special case of \cite[Theorem 6]{DM}, tells
us the ${\rm L}^2$-boundedness of a bounded $H_{\infty}$-calculus can be
extended to ${\rm L}^p$-boundedness, $p>1$.

\begin{lemma}\label{l2.2} Under the assumptions (a) and (b), the operator $L$ has a bounded
$H_\infty$-calculus in ${\rm L}^p({\mathbb R}^n)$, $p\in (1,\infty)$, that
is, $b(L)\in\mathcal{L}({\rm L}^p,{\rm L}^p)$ with
$$
\|b(L)\|_{p,p}=\|b(L)\|_{{\rm L}^p\to {\rm L}^p}\le c_{\nu,p}
\|b\|_\infty\quad\forall b\in H_{\infty}(S^0_{\nu}).
$$
Moreover, if $p=1$ then $b(L)$ is of weak type $(1,1)$.
\end{lemma}

\medskip
Thirdly, the Littlewood-Paley
 function   ${\mathcal G }_L(f)$ associated with an operator $L$
 is defined by
\begin{eqnarray}
{\mathcal G }_L(f)(x)=\Big( \int_0^{\infty}|\psi_t(L)f|^2\
{\frac{dt}{t}}\Big)^{1/2},
\label{e2.6}
\end{eqnarray}
where again $\psi\in \Psi(S^0_{\nu})$,
 and $\psi_t(\xi)=\psi(t\xi)$. It follows from Theorem 6 of \cite{AuDM}
 that the  function
${\mathcal G }_L(f)$ is bounded on   ${\rm L}^p$ for   $1<p<\infty$.
  More specifically,
 there exist constants $c_3,
c_4$  such that $0<c_3\leq c_4<\infty$ and
\begin{eqnarray}
c_3\|f\|_{{\rm L}^p}\leq \|{\mathcal G }_L(f)\|_{{\rm L}^p}\leq c_4\|f\|_{{\rm L}^p}
\label{e2.7}
\end{eqnarray}
for all  $f\in {\rm L}^p,  1<p<\infty$.

By duality, the operator  ${\mathcal G }_{L^{\ast}}(f)$ also satisfies
the estimate (\ref{e2.7}), where $L^{\ast}$ is the adjoint operator of $L$.

\medskip
\subsection{Acting class of semigroup $\{e^{-tL}\}_{t>0}$.}\ We now define the class
of functions that the operators $e^{-tL}$ act upon. Fix $1\leq p<\infty.$
For any  $\beta>0$, a complex-valued function $f\in {\rm L}^p_{\rm loc}({\mathbb R}^n)$ is
said to be a  function of type
$(p;  \beta)$   if $f$ satisfies
\begin{equation}
\left(\int_{{\mathbb R}^n}{\frac{|f(x)|^p}{(1+|x|)^{n+\beta}}}dx\right)^{1/p}\leq c<\infty.
\label{e2.8}
\end{equation}
We denote by ${\mathcal M}_{(p;  \beta)}$ the collection of all  functions
of type $(p; \beta)$. If $f\in {\mathcal M}_{(p;  \beta)},$  the norm
of $f\in {\mathcal M}_{(p;  \beta)}$ is defined by
$$
\|f\|_{{\mathcal M}_{(p;   \beta)}}
=\inf\{c\geq 0:\ \ (\ref{e2.8})  \ {\rm holds}\}.
$$
It is not hard to see that
  ${\mathcal M}_{(p;  \beta)}$ is a complex Banach space
under $
\|f\|_{{\mathcal M}_{(p;   \beta)}}<\infty$.
For any given  operator $L$, let
\begin{eqnarray}
\label{e2.9}
  {\Theta }(L)=\sup\big\{\epsilon>0:\ \ (\ref{e2.3})\ {\rm holds} \ \big\}
\end{eqnarray}
and write
 \begin{eqnarray*}
{\mathcal M}_p=
\left \{
\begin{array}{ll}
{\mathcal M}_{(p; {\Theta }(L))} &\ \ \ \ \ \  \ \ \ {\rm if } \ {\Theta }(L) <\infty;\\\\
  \bigcup\limits_{\beta:\
0<\beta<\infty}
{\mathcal M}_{(p;  \beta)} &\ \ \ \ \ \ \ \ \ {\rm if} \ {\Theta }(L)=\infty.
\end{array}
\right.
\end{eqnarray*}
Note that if $L=\triangle$ or $L=\sqrt{\triangle}$ on ${\mathbb R}^n$,
then ${\Theta }(\triangle)=\infty$ or ${\Theta }(\sqrt{\triangle})=1.$

For any $(x,t)\in {\mathbb R}^{n}\times (0, +\infty)={\mathbb R}^{n+1}_+$
and $f\in {\mathcal M}_p$, define
\begin{eqnarray}
P_tf(x)=e^{-tL}f(x)=\int_{{\mathbb R}^n} p_t(x,y)f(y)dy
\label{e2.10}
\end{eqnarray}
and
\begin{eqnarray}
\ \ \ Q_tf(x)=t Le^{-tL}f(x) =\int_{{\mathbb R}^n} -t\Big(\frac{d p_t(x,y)}{dt}\Big)f(y)dy.
\label{e2.11}
\end{eqnarray}
 It follows from  the estimate
(\ref{e2.4}) that
  the operators $P_tf$ and
$Q_tf$ are well defined. Moreover, the operator
$Q_t$ has the following two properties:

\smallskip

(i)\ For any $t_1, t_2>0$  and   almost all $x\in{\mathbb R}^n$,
$$Q_{t_1}Q_{t_2}f(x)=t_1t_2\Big(\frac{d^2 P_t}{dt^2}\big|_{t={t_1+t_2}}f\Big)(x);
$$

(ii)\    The kernel $q_{t^m}(x,y)$ of
$Q_{t^m}$ satisfies
\begin{eqnarray}
|q_{t^m}(x,y)|\leq c t^{-n}
g\Big( {|x-y|\over t}\Big)
\label{e2.12}\end{eqnarray}
where  the function $g$  satisfies the condition   (\ref{e2.3}).

\bigskip

\section{Basic properties}
\setcounter{equation}{0}

\subsection{A comparison of definitions} Assume that $L$ is an operator  which generates a semigroup
$e^{-tL}$ with the heat kernel bounds (\ref{e2.2}) and (\ref{e2.3}).
 In what follows,   $ B(x, {t})$ denotes the  ball
centered at $x$ and of the radius
${t}$.
Given  $B=B(x,  {t})$ and $\lambda>0$, we will write
$\lambda B$ for the $\lambda$-dilate
ball, which is the ball with the same center $x$ and
with radius $\lambda {t}$.

\begin{definition} \label{def3.1} Let $1\leq p<\infty$ and $\lambda\in (0, n).$ We say that

\item{\rm(i)} $f\in L_{loc}^p(\mathbb{R}^n)$ belongs to \LV\ provided (1.1) holds;

\item{\rm(ii)} $f\in {\mathcal M}_p$ associated with an operator $L$, is in
\LNV\ provided (1.2) holds.
\end{definition}

 \begin{remark}
\item{\rm (i)} Note first that (\LV, $\|\cdot \|_{{\rm L}^{p, \lambda}}$) and
(\LNV, $\|\cdot \|_{{\rm L}_{\rm L}^{p, \lambda}}$) are vector spaces with the
seminorms vanishing on constants and
$$
{\mathcal K}_{L,p}=\bigg\{f\in {\mathcal M}_p: P_tf(x)=f(x) \ {\rm for}\
 {\rm almost\
all\ }\ x\in {\mathbb R}^n {\rm \ for\ all}\ t>0 \bigg\},
$$
respectively. Of course, the spaces \LV and \LNV\  are understood to be modulo constants
and ${\mathcal K}_{L,p}$, respectively.  See Section 6 of \cite{DY2} for a discussion of
 the dimensions of ${\mathcal K}_{L,2}$ when $L$ is a second order elliptic operator of
 divergence form or a Schr\"odinger operator.

\item{\rm (ii)} We now give a list of examples of \LNV\ in different settings.

\ ($\alpha$) Define $P_t$ by putting $p_t(x,y)$ to be the heat kernel or the Poisson kernel:

$$
{1\over (4\pi t)^{n/2}}
e^{-|x-y|^2/4t}
\quad\hbox{or}\quad {c_n t\over (t^2 +|x-y|^2)^{(n+1)/2}}\quad\hbox{where}\quad c_n=
\frac{\Gamma\big(\frac{n+1}{2}\big)}{\pi^\frac{n+1}{2}}.
  $$
Then we will show that the corresponding space \LNV (modulo ${\mathcal
K}_{L,p}$) coincides with the classical \LV (modulo constants).

\smallskip

\ ($\beta$) Consider the Schr\"odinger operator with a non-negative potential $V(x)$:
$$
L=\triangle+V(x).
$$
To study singular integral operators associated to $L$ such as
functional calculi $f(L)$
or Riesz transform $\nabla L^{-1/2}$, it is useful to choose $P_t$
with kernel $p_t(x,y)$
to be the heat kernel (or Poisson kernel) of $L$.
By domination, its kernel $p_t(x,y)$ has a Gaussian upper bound (or a
Poisson bound).
\end{remark}

The following proposition shows that \LV\ is a subspace of \LNV\ in
many cases.

\begin{prop} \label{prop3.2}
Let $1\leq p<\infty$ and $\lambda\in (0, n).$
 Given an operator $L$ which generates a semigroup
$e^{-tL}$ with the heat kernel bounds (\ref{e2.2}) and (\ref{e2.3}).
A necessary and sufficient condition for the classical space
\LV\ $\subseteq$ \LNV\ with
\begin{equation}
\|f\|_{{\rm L}_L^{p, \lambda}}\leq c\|f\|_{{\rm L}^{p, \lambda}}
\label{e3.2}
\end{equation}
is that  for
every
$\
\! t>0$,
$e^{-tL}(1)=1$ almost everywhere, that is,
$
\int_{{\mathbb R}^n} p_t(x,y)dy=1$   for  almost  all
$x\in{\mathbb R}^n$.
\end{prop}
\begin{proof} Clearly, the condition $e^{-tL}(1)=1,$ a.e. is necessary for
\LV\ $\subseteq$ \LNV. Indeed, let us take
  $f=1$. Then, (\ref{e3.2}) implies
  $\|1\|_{{\rm L}_L^{p, \lambda}}=0$  and thus for every $t>0$,
$e^{-tL}(1)=1$ almost everywhere.

For the sufficiency, we borrow the idea of \cite[Proposition
3.1]{Ma}. To be more specific, suppose $f\in$ \LV\ . Then for any
Euclidean open ball $B$ with radius $r_B$, we compute
\begin{eqnarray*}
\|f-P_{t_B}f\|_{{\rm L}^p(B)}&\leq&\|f-f_B\|_{{\rm L}^p(B)}+\|f_B-P_{t_B}f\|_{{\rm L}^p(B)}\\
&\leq&\|f\|_{{\rm L}^{p,\lambda}}r_B^{\lambda/p}+\left(\int_B \Big(\int_{\mathbb R^n}
|f_B-f(y)|P_{t_B}(x,y)dy\Big)^pdx\right)^{1/p}\\
&=&\|f\|_{{\rm L}^{p,\lambda}}r_B^{{\lambda}/p}+ \left(\int_B \Big(I(B) +
J(B)\Big)^pdx\right)^{1/p},
\end{eqnarray*}
where
$$
I(B)=\int_{2B}|f_B-f(y)|P_{t_B}(x,y)dy
$$
and
$$
J(B)=\sum_{k=1}^\infty\int_{2^{k+1}B\setminus 2^k
B}|f_B-f(y)|P_{t_B}(x,y)dy.
$$
Next we make further estimates on $I(B)$ and $J(B)$. Thanks to
(\ref{e2.2}) and (\ref{e2.3}), we have
\begin{eqnarray*}
\|I(B)\|_{{\rm L}^p(B)}\le cr_B^{-n}g(0)\|f_B-f\|_{{\rm L}^1(B)}\leq c
r_B^{\lambda/p}\|f\|_{{\rm L}^{p,\lambda}}.
\end{eqnarray*}
Again, using (\ref{e2.2}) and (\ref{e2.3}), we derive that for $x\in
B$ and $y\in 2^{k+1}B\setminus 2^kB$,
$$
P_{t_B}(x,y)\le c r_B^{-n}g(2^k)\le
cr_B^{-n}2^{-k(n+\epsilon)},\quad k=1,2,....
$$
where $\epsilon>0$ is a constant. Consequently,
\begin{eqnarray*}
\|J(B)\|_{{\rm L}^p(B)}&\le& c r_B^{-n}\left(\int_B\left(\sum_{k=1}^\infty
g(2^k)\int_{2^{k+1}B\setminus
2^kB}|f_B-f(y)|dy\right)^pdx\right)^{1/p}\\
&\le& cr_B^{n/p-n}\sum_{k=1}^\infty
g(2^k)\left(\int_{2^{k+1}B}|f_{2^{k+1}B}-f(y)|dy+
(2^kr_B)^n|f_{2^{k+1}B}-f_B|\right)\\
&\le&c r_B^{\lambda/p}\|f\|_{{\rm L}^{p,\lambda}}\left(\sum_{k=1}^\infty
2^{-k(\epsilon+\frac{n-\lambda}{p})}+\sum_{k=1}^\infty
k2^{-k\epsilon}\right).
\end{eqnarray*}
Putting these inequalities together, we find $f\in$\ \LNV.
\end{proof}

\subsection{Fundamental characterizations}

In the argument for Proposition \ref{prop3.2}, we have used the
following crucial fact that for any $f\in$\LV\ and a constant $K>1$,
\begin{eqnarray*}
|f_B-f_{KB}|\leq cr_B^{\lambda-n\over p}\|f\|_{{\rm L}^{p,
\lambda}}.
\end{eqnarray*}

Now, this property can be used to give a characterization of \LV\
spaces in terms of the Poisson integral. To this end, we denote the
Laplacian by $\triangle=-\sum_{i=1}^n\partial^2_{x_i}$
  and  $e^{-t\sqrt{\triangle}}$  the Poisson semigroup
  on ${\mathbb R}^n$. We observe that if
$$
f\in {\mathcal M}_{\sqrt{\triangle}, p}=\bigg\{f\in {\rm L}^p_{\rm
loc}({\mathbb R}^n): |f(x)|^p(1+|x|^{n+1})^{-1}\in {\rm L}^1({\mathbb
R}^n)\bigg\},
$$
then we can define the operator $e^{-t\sqrt{\triangle}}$ by the
  Poisson integral as follows:
$$
 e^{-t\sqrt{\triangle}}f(x)=\int_{{\mathbb R}^n}p_t(x-y)f(y)dy,\ \ \ t>0,
$$
where
$$
p_t(x-y)={c_n t\over (t^2+|x-y|^2)^{(n+1)/2}}.
$$

\begin{prop}\label{prop3.2a} Let $1\leq p<\infty$, $\lambda\in (0,n)$ and
$f\in {\mathcal M}_{\sqrt{\triangle}, p}$. Then
$f\in$ \LV\ if and only if
\begin{eqnarray}
|\!|\!|f|\!|\!|_{{\rm L}^{p,\lambda}({\mathbb
R}^n)}=\Big(\sup_{(x,t)\in{\mathbb
R}^{n+1}_+}t^{n-\lambda}e^{-t\sqrt{\Delta}}\big(|f-e^{-t\sqrt{\Delta}}f(x)|^p\big)(x)\Big)^{1/p}<\infty.
\label{e3.2bc}
\end{eqnarray}
\end{prop}
\begin{proof} On the one hand, assume (\ref{e3.2bc}). Note that
$|y-x|<t$ implies
$$
\frac{c_n t}{(t^2+|y-x|^2)^{\frac{n+1}{2}}}\ge c t^{-n}.
$$
For a fixed ball $B=B(x, r_B)$ centered at $x$, we let $t_B=r_B$. We then  have
\begin{eqnarray*}
r_B^{-\lambda}\|f-f_B\|_{{\rm L}^p(B)}^p&\leq&
cr_B^{-\lambda}\|f-e^{-t_B\sqrt{\Delta}}f(x)\|_{{\rm L}^p(B)}^p\\
&\leq& c r_B^{n-\lambda}\int_{B}|f(y)-e^{-t_B\sqrt{\Delta}}f(x)|^p
\frac{ c_nt_B}{(t_B^2+|y-x|^2)^{\frac{n+1}{2}}}dy\\
&\leq& c|\!|\!|f|\!|\!|_{{\rm L}^{p,\lambda}}^p,
\end{eqnarray*}
whence producing $f\in$ \LV.

On the other hand, suppose $f\in$ \LV. \ In a similar manner to
proving the sufficiency part of Proposition \ref{prop3.2}, we obtain
that if $(x,t)\in{\mathbb R}^{n+1}_+$ then
\begin{eqnarray*}
e^{-t\sqrt{\Delta}}\big(|f-e^{-t\sqrt{\Delta}}f(x)|^p\big)(x)&\le& c
t^{\lambda-n}\|f\|_{{\rm L}^{p,\lambda}}^p+c\sum_{j=1}^\infty\int_{2^{k+1}B\setminus
2^kB}\frac{|f(y)-f_B|^p t}{(t^2+|y-x|^2)^{\frac{n+1}{2}}}dy\\
&\leq& c t^{\lambda-n}\|f\|_{{\rm L}^{p,\lambda}}^p,
\end{eqnarray*}
and hence (\ref{e3.2bc}) holds.
\end{proof}

\begin{remark} Since a simple computation gives
\begin{eqnarray*}
&&e^{-t\sqrt{\Delta}}\big(|f-e^{-t\sqrt{\Delta}}f(x)|^2\big)(x)\\
&&=\int_{{\mathbb
R}^n}\big(f(y)-e^{-t\sqrt{\Delta}}f(x)\big)\overline{\big(f(y)-e^{-t\sqrt{\Delta}}f(x)\big)}p_t(x-y)dy\\
&&=\int_{{\mathbb
R}^n}|f(y)|^2p_t(x-y)dy-e^{-t\sqrt{\Delta}}f(x)\left(\int_{{\mathbb
R}^n}\overline{f(y)}p_t(x-y)dy\right)\\
&&\quad -\overline{e^{-t\sqrt{\Delta}}f(x)}\left(\int_{{\mathbb
R}^n}{f(y)}p_t(x-y)dy\right)+\big|e^{-t\sqrt{\Delta}}f(x)\big|^2\\
&&=
e^{-t\sqrt{\Delta}}|f|^2(x)-\big|e^{-t\sqrt{\Delta}}f(x)\big|^2,
\end{eqnarray*}
we have that if $f\in {\mathcal M}_{\sqrt{\triangle}, 2}$ then
$f\in {\rm L}^{2,\lambda}({\mathbb R}^{n})$ when and only when
$$
\sup_{(x,t)\in{\mathbb
R}^{n+1}_+}t^{n-\lambda}\Big(e^{-t\sqrt{\Delta}}|f|^2(x)-\big|e^{-t\sqrt{\Delta}}f(x)\big|^2\Big)<\infty
$$
which is equivalent to (see also \cite{Lu} for the BMO-setting, i.e., $\lambda=n$)
$$
\sup_{(x,t)\in{\mathbb
R}^{n+1}_+}t^{n-\lambda}\int_{{\mathbb
R}^{n+1}_+}G_{{\mathbb
R}^{n+1}_+}\big((x,t), (y,s)\big)|\nabla_{y,s} e^{-s\sqrt{\Delta}}f(y)|^2dyds<\infty,
$$
where $G_{{\mathbb R}^{n+1}_+}\big((x,t), (y,s)\big)$ is the Green function of
${\mathbb R}^{n+1}_+$ and $\nabla_{y,s}$ is the gradient operator in the space-time variable $(y,s)$.
\end{remark}

To find out an\ \LNV\ analog of Proposition \ref{prop3.2a}, we take
 Proposition 2.6 of \cite{DY1} into account, and establish
 the following  property of the  class of operators $P_t$.

\begin{lemma} \label{e3.3} Let $1\leq p<\infty$ and $\lambda\in (0, n)$.
Suppose $f\in$ \LNV.  Then:
\item{\rm(i)} For any $t>0$ and $K>1$, there exists  a constant  $c>0$
independent of
 $t$ and $K$ such that
\begin{eqnarray}
|P_tf(x)-P_{Kt} f(x)|\leq c t^{\lambda-n\over pm}\|f\|_{{\rm
L}_L^{p, \lambda}} \label{e3.3n}
\end{eqnarray}
for almost all $x\in {\mathbb R}^n$.
\item{\rm(ii)} For any $\delta>0$, there exists $c(\delta)>0$ such that
\begin{eqnarray}
\int_{{\mathbb R}^n} {t^{\delta/m}\over
(t^{1/m}+|x-y)^{n+\delta}}|({\mathcal I}-P_t)f(y)|dy \leq {c(\delta)
t^{\frac{\lambda-n}{pm}}}\|f\|_{{\rm L}_L^{p, \lambda}}\label{e3.31}
\end{eqnarray}
for any $x\in {\mathbb R}^n$,
\end{lemma}

\begin{proof} (i) For any $t>0$, we choose $s$ such that $t/4\leq
s\leq t$.
Assume that  $f\in$ \LNV, where  $1\leq p<\infty$ and $\lambda\in (0, n)$, we
estimate the term $|P_tf(x)-P_{t+s} f(x)|$. Using the
commutative property of the semigroup $\{P_t\}_{t>0}$,
 we can write
$$
P_tf(x)-P_{t+s} f(x)=P_t (f-P_{s} f)(x).
$$
Since $f\in$ \LNV, one has
\begin{eqnarray*}
|P_tf(x)-P_{t+s}f(x)|&\leq&\int_{{\mathbb R}^n}|p_t(x,y)||f(y)-P_{s}f(y)| dy\\
&\leq&{c\over |B(x, t^{1/m})|}\int_{{\mathbb R}^n}
\Big(1+{|x-y|\over t^{1/m}}\Big)^{-(n+\epsilon)}|f(y)-P_{s}f(y)| dy\\
&\leq&c\Big({1\over |B(x, s^{1/m})|}\int_{B(x, s^{1/m})}|f(y)-P_sf(y)|^p dy\Big)^{1/p}\\
&&+{c\over |B(x, s^{1/m})|}\int_{B(x,
s^{1/m})^c}\Big(1+{|x-y|\over s^{1/m}}\Big)^{-(n+\epsilon)}|f(y)-P_sf(y)| dy\\
&\leq&cs^{\lambda-n\over pm}\|f\|_{{\rm L}_L^{p, \lambda}} +{\rm I}.
\end{eqnarray*}
We then decompose ${\mathbb R}^n$ into a geometrically increasing sequence
of concentric balls, and obtain
\begin{eqnarray*}
{\rm I}&=&
c\sum_{k=0}^{\infty}{1\over |B(x, s^{1/m})|}\int_{B(x, 2^{k+1} s^{1/m})\backslash B(x, 2^k s^{1/m})}
{\Big(1+{|x-y|\over s^{1/m}}\Big)^{-(n+\epsilon)}} |f(y)-P_sf(y)| dy\\
&\leq&c\sum_{k=0}^{\infty}2^{-k(n+\epsilon)}{1\over |B(x,  s^{1/m})|}\int_{B(x, 2^{k+1}
s^{1/m})}
|f(y)-P_s f(y)| dy,
\end{eqnarray*}
since
$$
{\big(1+{s^{-1/m}|x-y|}\big)^{-n-\epsilon}}\leq c2^{-k(n+\epsilon)}\quad\forall\ y\in
B(x, 2^{k+1} s^{1/m})\backslash
B(x, 2^k s^{1/m}).
$$

For a fixed
positive integer $k$, we consider the ball $B(x, 2^ks^{1/m})$. This ball is contained in the
cube $Q[x, 2^{k+1}s^{1/m}]$  centered at $x$ and of the side length
$2^{k+1}s^{1/m}$.
We then divide this cube $Q[x, 2^{k+1}s^{1/m}]$ into    $[2^{k+1}([\sqrt{n}\
\!]+1)]^n$ small cubes $\{Q_{x_{k_i}}\}_{i=1}^{N_k}$ centered at $x_{k_i}$ and of equal
side length
$([\sqrt{n}\
\!]+1)^{-1}s^{1/m}$, where $N_k=[2^{k+1}([\sqrt{n}\ \!]+1)]^n$. For any  $i=1, 2,
\cdots,  N_k$, each of these small cubes $Q_{x_{k_i}}$   is
then contained in the corresponding  ball
$B_{k_i}$ with the same center ${x_{k_i}}$ and radius
$r=s^{1/m}$,   Consequently, for any ball $B(x, 2^kt)$,  $k=1, 2, \cdots,
$ there exists  a corresponding collection of balls  $B_{k_1}, B_{k_2},
\cdots, B_{k_{N_k}}$ such that

\smallskip

(i)\ each ball $B_{k_i}$ is of the radius t;

\smallskip

(ii)\ $B(x, 2^{k} s^{1/m})\subset \bigcup\limits_{i=1}^{N_k}B_{k_i};$

\smallskip

(iii)\ there exists a constant $c>0$ independent of $k$ such that $N_k\leq c2^{kn};$
\smallskip

(iv)\ each point of $B(x, 2^{k} s^{1/m})$ is contained in at most
a finite number $M$ of the balls $B_{k_i}$, where $M$ is independent of $k$.

\smallskip

Applying  the properties (i), (ii), (iii) and (iv)
above, we obtain
\begin{eqnarray*}
{\rm I}
&\leq&c\sum_{k=0}^{\infty}2^{-k(n+\epsilon)}{1\over |B(x,  s^{1/m})|}
\int_{\bigcup\limits_{i=1}^{N_{k+1}}B_{k_i}}
|f(y)-P_t f(y)| dy\\
&\leq&c\sum_{k=0}^{\infty}2^{-k(n+\epsilon)}\sum_{i=1}^{N_{k+1}}
{1\over |B_{k_i}|}\int_{B_{k_i}}
|f(y)-P_s f(y)| dy\\
&\leq&c\sum_{k=0}^{\infty}2^{-k(n+\epsilon)}N_{k+1}\sup_{i: 1\leq i\leq N_{k+1}}
\bigg({1\over |B_{k_i}|}\int_{B_{k_i}}
|f(y)-P_s f(y)|^p dy\bigg)^{1/p}\\
&\leq&c\sum_{k=0}^{\infty}2^{-k(n+\epsilon)}2^{kn}s^{\lambda-n\over pm}
\|f\|_{{\rm L}_L^{p, \lambda}}\\
&\leq&cs^{\lambda-n\over pm}\|f\|_{{\rm L}_L^{p, \lambda}},
\end{eqnarray*}
which gives (\ref{e3.3n}) for the case $t/4\leq s\leq t.$

For the case $0<s<t/4$, we write
$$P_t f(x)-P_{t+s}f(x)=(P_t  f(x)-P_{2t}  f(x))
-(P_{t+s} (f-P_{t-s}f)(x).
$$
Noting that $(t+s)/4\leq (t-s)<t+s$, we obtain (\ref{e3.3n}) by using the
same argument as above. In general, for
any $K>1$, let $l$ be the integer satisfying $ 2^l\leq K<2^{l+1},
$ hence $l\leq 2{\rm log} K$.  This, together with the fact that
$\lambda\in (0, n)$, imply that there exists a constant  $c>0$
independent of $t$ and $K$ such that
\begin{eqnarray*}
|P_t
f(x)-P_{Kt}  f(x)|&\leq&\sum_{k=0}^{l-1}|P_{2^kt}
f(x)-P_{2^{k+1}t}  f(x)| +|P_{2^{l}t} f(x)-P_{Kt} f(x)|\\
&\leq&c\sum_{k=0}^{l-1}(2^kt)^{\lambda-n\over pm}\|f\|_{{\rm L}_L^{p, \lambda}}
+c(Kt)^{\lambda-n\over pm}\|f\|_{{\rm L}_L^{p, \lambda}}\\
&\leq&
ct^{\lambda-n\over pm} \|f\|_{{\rm L}_L^{p, \lambda}}
\end{eqnarray*}
for almost all $x\in {\mathbb R}^n$.

(ii) Choosing a ball $B$ centered at $x$ and of the radius
$r_B=t^{1/m}$, and using (\ref{e3.3n}), we have
   \begin{eqnarray}
   \Big({1\over |2^kB|}&&\hspace{-1cm} \int_{2^{k}B}|f(y)-P_{t}f(y)|^p
dy\Big)^{1/p}\nonumber\\
   &\leq&\Big({1\over |2^kB|}\int_{2^{k}B}|f(y)-P_{t_{2^kB}}f(y)|^p
dy\Big)^{1/p}
   +\sup_{y\in 2^kB}|P_{t_{2^kB}}f(y)-P_{t}f(y)| \nonumber\\
   &\leq& ct^{\lambda-n\over pm} \|f\|_{{\rm L}_L^{p, \lambda}}
   \label{e3.4}
   \end{eqnarray}
   for all $k$. Putting $2^{-1}B=\emptyset$, we read off
   \begin{eqnarray*}
\int_{{\mathbb R}^n} &&\hspace{-1cm}{t^{\delta/m}\over
(t^{1/m}+|x-y)^{n+\delta}}|({\mathcal I}-P_t)f(y)|dy\\
   &\leq&\sum_{k=0}^{\infty} \int_{2^kB\backslash 2^{k-1}B}
{t^{\delta/m}\over
(t^{1/m}+|x-y)^{n+\delta}}|({\mathcal I}-P_t)f(y)|dy\nonumber\\
&\leq&c\sum_{k=0}^{\infty} 2^{kn}2^{-k(n+\delta)}
{1\over |2^kB|}\int_{2^{k}B} |f(y)-P_{t}f(y)|dy\nonumber\\
&\leq&c\sum_{k=0}^{\infty} 2^{-k \delta}
\Big({1\over |2^kB|}\int_{2^{k}B} |f(y)-P_{t}f(y)|^pdy\Big)^{1/p}
\nonumber\\
&\leq& c\sum_{k=0}^{\infty} 2^{-k\delta}
t^{\frac{\lambda-n}{pm}}\|f\|_{{\rm L}_L^{p, \lambda}}\nonumber\\
&\leq&ct^{\frac{\lambda-n}{pm}}\|f\|_{{\rm L}_L^{p, \lambda}}.
\end{eqnarray*}
\end{proof}

The above analysis suggests us to introduce the maximal Morrey space as follows.

\begin{definition}\label{definition 3.5}  Let $1\leq p<\infty$
and  $\lambda\in (0, n)$.
We say that $f\in {\mathcal M}_p$ is in
${\rm L}^{p,\lambda}_{L, {\rm max}}({\mathbb R}^n)$
associated with an operator $L$,
if there exists some constant $c$ (depending on $f$) such that
\begin{equation}
  \big|P_t(|f-P_{t}f|^p)(x)\big|^{1/p}\leq ct^{\lambda-n\over pm}\quad\hbox{for\
  almost\ all}\ x\in{\mathbb R}^n\ \ \hbox{and}\ \ t>0.
\label{e3.6a}
\end{equation}
The smallest bound $c$
for which (\ref{e3.6a}) holds then taken to be the norm of
$f$ in this space, and is denoted by $\|f\|_{{\rm L}^{p,\lambda}_{L,
{\rm max}}}$.
\end{definition}

Using Lemma \ref{e3.3}, we can derive a characterization in terms of the maximal
 Morrey space under an extra hypothesis.

\begin{prop}\label{prop 3.6} Let $1\leq p<\infty$ and  $\lambda\in (0, n)$. Given
an operator $L$ which generates a semigroup $e^{-tL}$ with the heat kernel
bounds (2.2) and (2.3). Then \LNV\ $\subseteq{\rm L}^{p,\lambda}_{L,
{\rm max}}({\mathbb R}^n)$. Furthermore, if the kernels $p_t(x,y)$
of operators $P_t$ are nonnegative functions when $t>0$, and satisfy
the following lower bounds
\begin{equation}
p_t(x,y)\geq {c\over t^{n/m}}
\label{e3.6}
\end{equation}
for   some positive constant $c$independent of $t$,    $x$ and $y$,
then, ${\rm L}^{p,\lambda}_{L, {\rm max}}({\mathbb R}^n)=$ \LNV\ .
\end{prop}

\begin{proof} Let us first prove \LNV\ $\subseteq$
${\rm L}^{p,\lambda}_{L, {\rm max}}({\mathbb R}^n).$
For any fixed
$t>0$ and $x\in {\mathbb R}^n$,
we choose a ball $B$ centered at $x$ and of the radius
$r_B=t^{1/m}$.  Let $f\in$ \LNV.
To estimate (\ref{e3.6a}), we use the decay of function
$g$ in (\ref{e2.3}) to get
   \begin{eqnarray*}
\big|P_t(|f-P_{t}f|^p)(x)\big| &\leq&\int_{{\mathbb R}^n}|p_t(x,y)||f(y)-P_{t}f(y)|^pdy
\nonumber\\
   &\leq&c\sum_{k=0}^{\infty} {1\over |B|}\int_{2^kB\backslash 2^{k-1}B}
g\bigg({|x-y|\over t^{1/m}}\bigg) |f(y)-P_{t}f(y)|^pdy\nonumber\\
&\leq&c\sum_{k=0}^{\infty} 2^{kn}g(2^{(k-1) })
{1\over |2^kB|}\int_{2^{k}B} |f(y)-P_{t}f(y)|^pdy\nonumber\\
&\leq& c\sum_{k=0}^{\infty}2^{kn}g(2^{(k-1) })
t^{\lambda-n\over  m} \|f\|^p_{{\rm L}_L^{p, \lambda}}\nonumber\\
&\leq&ct^{\lambda-n\over  m} \|f\|^p_{{\rm L}_L^{p, \lambda}}.
\end{eqnarray*}
 This proves $\|f\|_{{\rm L}_{L, {\rm max}}^{p, \lambda}}
\leq c\|f\|_{{\rm L}_L^{p, \lambda}}$.

We now prove ${\rm L}^{p,\lambda}_{L, {\rm max}}({\mathbb R}^n)$
 $\subseteq$ \LNV\ under (\ref{e3.6}). For a fixed ball $B=B(x, r_B)$
centered at $x$, we
let $t_B=r^m_B$. For any $f\in {\rm L}^{p,
\lambda}_{L, {\rm max}}({\mathbb R}^n)$, it follows from (\ref{e3.6}) that
one has
\begin{eqnarray*}
   {1\over |B|} \int_B|f(y)-P_{t_B}f(y)|^pdy
&\leq& c\int_{B(x, t_B^{1/m})}p_{t_B}(x,y)|f(y)-P_{t_B}f(y)|^pdy
\\
&\leq& c\int_{{\mathbb R}^n}p_{t_B}(x,y)|f(y)-P_{t_B}f(y)|^pdy\\
&\leq&ct_B^{\lambda-n\over m}\|f\|^p_{{\rm L}_{L, {\rm max}}^{p, \lambda}},
\end{eqnarray*}
which proves $\|f\|_{{\rm L}_L^{p, \lambda}}\leq
c\|f\|_{{\rm L}_{L, {\rm max}}^{p, \lambda}}
$. Hence,
the proof of Proposition \ref{prop 3.6} is complete.
\end{proof}

\section{An identity for the dual pairing}
\setcounter{equation}{0}

\subsection{A dual inequality and a reproducing formula} From now on, we need the following notation.
Suppose $B$ is an open ball centered at $x_B$ with radius $r_B$ and $f\in {\mathcal M}_{p}$.
Given an ${\rm L}^q$ function $g$ supported on a ball $B$, where ${1\over q}+ {1\over p}=1$. For
 any $ (x,t)\in {\mathbb R}^{n+1}_+$, let
\begin{eqnarray}
\hspace{1.5cm} F(x,t)=Q_{t^m} ({\mathcal {I}}-P_{t^m})f(x)\ \ \ {\rm and}\ \ \
G(x,t)=Q^{\ast}_{t^m} ({\mathcal {I}}-P^{\ast}_{r_B^m})g(x),
\label{e4.6}
\end{eqnarray}
where $P^{\ast}_t$ and $Q^{\ast}_t$  are the adjoint operators of $P_t$
and $Q_t$, respectively.

\begin{lemma} \label{lemma 4.4} Assume that $L$ satisfies the
assumptions (a) and (b) of Section 2.2. Suppose $f, g, F, G, p, q$ are as in (\ref{e4.6}).

\item{\rm(i)} If $f$ also satisfies
$$
|\!|\!|f|\!|\!|_{{\rm L}_L^{p,\lambda}}=\sup_{B\subset{\mathbb R}^n}
r_B^{-\frac{\lambda}{p}}\Big\|\Big\{\int_0^{r_B}|Q_{t^m}(I-P_{t^m})f(x)|^2
\frac{dt}{t}\Big\}^{1/2}\Big\|_{{\rm L}^p(B)}<\infty,
$$
where the supremum is taken over all open ball $B\subset{\mathbb R}^n$ with radius
$r_B$, then there exists a constant $c>0$ independent of any open ball $B$ with radius $r_B$ such that
\begin{eqnarray}
\int_{{\mathbb R}^{n+1}_+}
|F(x,t) G(x,t)|{dxdt\over t}\leq c r_B^{\lambda/p}|\!|\!|f|\!|\!|_{{\rm L}_L^{p,\lambda}}
\|g\|_{{\rm L}^q}.
\label{e4.7}
\end{eqnarray}

\item{\rm(ii)} If
$$
h\in {\rm L}^q({\mathbb R}^n),\quad b_m=\frac{36m}{5}\quad\hbox{and}
\quad
1=b_m\int_0^{\infty}t^{2m}e^{-2t^m}(1-e^{-t^m}){dt\over t},
$$
then
$$
h(x)=b_m\int_0^{\infty}
(Q^{\ast}_{t^m})^2 ({\mathcal {I}}-P^{\ast}_{t^m})h(x){dt\over t},
$$
where the integral converges strongly in ${\rm L}^q({\mathbb R}^n)$.
\end{lemma}

\begin{proof} (i) For any ball $B\subset{\mathbb R}^n$ with radius $r_B$, we still put
$$
T(B)=\{(x,t)\in{\mathbb R}^{n+1}_+: x\in B, \ 0<t<r_B\}.
$$
We then write
\begin{eqnarray*}
\int_{{\mathbb R}^{n+1}_+}
|F(x,t) G(x,t)|{dxdt\over t} &=&\int_{T(4B)} \big|F(x,t)
G(x,t)\big|{dxdt\over t}\\
&&+\sum_{k=1}^{\infty}\int_{T(2^{k+1}B)\backslash T(2^{k}B) }
\big|F(x,t) G(x,t)\big|{dxdt\over t}\\
&=&{\rm A_1} + \sum_{k=2}^{\infty} {\rm A_k}
\end{eqnarray*}
Recall that $q>1$ and
${1\over q}+ {1\over p}=1$. Using the H\"older inequality, together with
(\ref{e2.7}) (here
$\psi(z)=ze^{-z}$), we obtain
\begin{eqnarray*}
{\rm A_1}&\leq&\Big\|\Big\{\int_0^{r_{2B}}\big|Q_{t^m} ({\mathcal
{I}}-P_{t^m})f(x) \big|^2{dt\over t}\Big{\}}^{1/2}\Big\|_{{\rm L}^p(2B)}\\
&&\times\Big\|\Big\{\int_0^{r_{2B}}\big|Q_{t^m}^{\ast}({\mathcal
{I}}-P^{\ast}_{r_B^m})g(x) \big|^2{dt\over
t}\Big{\}}^{1/2}\Big\|_{{\rm L}^q(2B)}
\\
&\leq&\|\Big\{\int_0^{r_{2B}}\big|Q_{t^m} ({\mathcal {I}}-P_{t^m})
f(x)\big|^2{dt\over t}\Big{\}}^{1/2}\|_{{\rm L}^p(2B)}
\|{\mathcal G}_{L^{\ast}}(({\mathcal I}-P^{\ast}_{r_B^m})g)\|_{{\rm L}^q}\\
&\leq&cr_B^{\lambda \over p}|\!|\!|f|\!|\!|_{{\rm L}_L^{p,\lambda}}\|g\|_{{\rm L}^q}.
\end{eqnarray*}
Let us estimate ${\rm A_k}$ for $k=2,3, \cdots.$
Note that for   $x\in T(2^{k+1}B)\backslash T(2^{k}B)$
and $y\in B$, we have that $|x-y|\geq 2^kr_B$. Using (\ref{e2.4}) and the commutative property
of $\{P_t\}_{t>0}$, we get
 \begin{eqnarray*}
|Q_{t^m}^{\ast} ({\mathcal {I}}-P^{\ast}_{r_B^m})g(x)|&\leq & |Q_{t^m}^{\ast}g(x)|
+c\Big({t\over t+r_B}\Big)^m\Big|Q_{t^m+r_B^m}g(x)\Big|\\
&\leq&c\int_{B} {t^{\epsilon}
\over (t+ |x-y|)^{n+{\epsilon}}}|g(y)|dy\\
&&+c\Big({t\over r_B}\Big)^m\int_{B} {r_B^{\epsilon}
\over (r_B+ |x-y|)^{n+{\epsilon}}}|g(y)|dy\\
 &\leq& c  {t^{\epsilon_0}
\over (2^kr_B)^{n+{\epsilon_0} }}\int_{B}|g(y)|dy\\
&\leq& c  {t^{\epsilon_0}
\over (2^kr_B)^{n+ {\epsilon_0}}}r_B^{n\over p}\|g\|_{{\rm L}^q},
\end{eqnarray*}
where ${\epsilon_0}=2^{-1}{\rm min} (m, \epsilon)$ and $q=p/(p-1)$.
Consequently,
\begin{eqnarray*}
\Big{\|}\Big\{\int_0^{2^kr_B}
|Q_{t^m}^{\ast}({\mathcal {I}}-P^{\ast}_{r_B^m})g(x)\chi_{T(2^{k+1}B)\backslash T(2^{k}B) }|^2
{dt\over t}\Big\}^{1/2}\Big{\|}_{{\rm L}^{q}(2^kB)}
\leq c2^{kn( {1\over q}-1)}\|g\|_{{\rm L}^q}.
\end{eqnarray*}
Therefore,
\begin{eqnarray*}
{\rm A}_k
&\leq&\Big{\|}\Big\{\int_0^{2^kr_B}|Q_{t^m} ({\mathcal {I}}-P_{t^m}) f(x)
|^2{dt\over t}\Big\}^{1/2}\Big{\|}_{{\rm L}^p(2^kB)}\\
&&\times \Big{\|}
\Big\{\int_0^{2^kr_B}|Q_{t^m}^{\ast}({\mathcal {I}}-P^{\ast}_{r_B^m})
g(x)\chi_{T(2^{k+1}B)\backslash T(2^{k}B) }|^2
{dt\over t}\Big\}^{1/2}\Big{\|}_{{\rm L}^{q}(2^kB)}\\
&\leq&c (2^kr_B)^{\lambda \over p}2^{kn({1\over q}-1)}|\!|\!|f|\!|\!|_{{\rm L}_L^{p,\lambda}}
\|g\|_{{\rm L}^q}\\
&\leq&c2^{k(\lambda- n)\over p}r_B^{{\lambda  \over p}}|\!|\!|f|\!|\!|_{{\rm L}_L^{p,\lambda}}
\|g\|_{{\rm L}^q}.
\end{eqnarray*}
Since $\lambda\in (0, n)$, we  have
\begin{eqnarray*}
\int_{{\mathbb R}^{n+1}_+}
|F(x,t) G(x,t)|{dxdt\over t}
&\leq& cr_B^{{\lambda  \over p}}|\!|\!|f|\!|\!|_{{\rm L}_L^{p,\lambda}}
\|g\|_{{\rm L}^q}+
c\sum_{k=1}^{\infty} 2^{k(\lambda-n)\over 2} r_B^{{\lambda  \over p}}
|\!|\!|f|\!|\!|_{{\rm L}_L^{p,\lambda}}\|g\|_{{\rm L}^q}\\
&\leq& cr_B^{{\lambda  \over p}}
|\!|\!|f|\!|\!|_{{\rm L}_L^{p,\lambda}}\|g\|_{{\rm L}^q},
\end{eqnarray*}
as desired.

(ii) From Lemma \ref{l2.2} we know that $L$ has a bounded $H_{\infty}$-calculus in ${\rm L}^q$ for all $q>1$.
This, together with elementary integration, shows that $\{g_{\alpha
\beta}(L^{\ast})\}$ is a uniformly bounded net in ${\mathcal L}({\rm L}^q,
{\rm L}^q)$, where
$$
g_{\alpha \beta}(L^{\ast})=b_m\int_{\alpha}^{\beta}
(Q^{\ast}_{t^m})^2 ({\mathcal {I}}-P^{\ast}_{t^m}){dt\over t}
$$
for all $0<\alpha<\beta<\infty.$

As a consequence of Lemma \ref{l2.1}, we have that for any $h\in
{\rm L}^q({\mathbb R}^n)$,
$$
h(x)=b_m\int_0^{\infty}
(Q^{\ast}_{t^m})^2 ({\mathcal {I}}-P^{\ast}_{t^m})h(x){dt\over t}
$$
where $b_m=\frac{36m}{5}$ and the integral is strongly convergent in ${\rm L}^q({\mathbb R}^n)$.
\end{proof}

\subsection{The desired dual identity} Next, we establish the following dual identity
 associated with the operator $L$.

\begin{prop} \label{prop 4.3}  Assume that $L$ satisfies the
assumptions  (a) and  (b) of Section 2.2. Suppose $B, f, g, F, G, p, q$ are
defined as in (\ref{e4.6}). If $|\!|\!|f|\!|\!|_{{\rm L}_L^{p,\lambda}}<\infty$ and
$b_m=\frac{36m}{5}$, then
\begin{eqnarray}
\int_{{\mathbb R}^n} f(x) ({\mathcal {I}}-P^{\ast}_{r_B^m})g(x)dx=b_m\int_{{\mathbb R}^{n+1}_+}
F(x,t) G(x,t){dxdt\over t}.
\label{e4.8}
\end{eqnarray}
\end{prop}
\begin{proof} From Lemma \ref{lemma 4.4} (i) it turns out that
 $$
\int_{{\mathbb R}^{n+1}_+} \big|F(x,t) G(x,t)\big| {dxdt/t} <\infty.
$$
By the dominated convergence theorem, the following integral
converges  absolutely and satisfies
$$
\int_{{\mathbb R}^{n+1}_+} F(x,t) G(x,t)
{dxdt\over t}=\lim_{\delta\rightarrow 0}\lim_{N\rightarrow \infty}
\int_{\delta}^N\int_{{\mathbb R}^n} F(x,t) G(x,t)
{dxdt\over t}.
$$
Next, by Fubini's theorem, together with
 the commutative property of the semigroup $\{e^{-tL}\}_{t>0}$,
we have
$$
\int_{{\mathbb R}^n} F(x,t) G(x,t)dx=\int_{{\mathbb R}^{n}} f(x) (Q^{\ast}_{t^m})^2
({\mathcal {I}}-P^{\ast}_{t^m})({\mathcal {I}}-P^{\ast}_{r_B^m})g(x)dx,
\ \ \ \ \ \forall t>0.
$$
This gives
\begin{eqnarray}
\label{e4.9}
 \int_{{\mathbb R}^{n+1}_+}&&\hspace{-1cm} F(x,t) G(x,t)
{dxdt\over t}\nonumber\\
&=&\lim_{\delta\rightarrow 0}\lim_{N\rightarrow \infty}
 \int_{\delta}^N\Big[\int_{{\mathbb R}^n} f(x)(Q^{\ast}_{t^m})^2
({\mathcal {I}}-P^{\ast}_{t^m})({\mathcal {I}}-P^{\ast}_{r_B^m})
g(x)dx\Big]{dt\over t}\nonumber\\
&=&\lim_{\delta\rightarrow 0}\lim_{N\rightarrow \infty}
 \int_{{\mathbb R}^n} f(x)\Big[\int_{\delta}^N (Q^{\ast}_{t^m})^2 ({\mathcal
{I}}-P^{\ast}_{t^m})({\mathcal {I}}-P^{\ast}_{r_B^m})g(x){dt\over t}\Big]dx\nonumber\\
&=&\lim_{\delta\rightarrow 0}\lim_{N\rightarrow \infty}
 \int_{{\mathbb R}^n} f_1(x)\Big[\int_{\delta}^N (Q^{\ast}_{t^m})^2 ({\mathcal
{I}}-P^{\ast}_{t^m})({\mathcal {I}}-P^{\ast}_{r_B^m})g(x){dt\over t}\Big]dx \nonumber\\
&&+ \lim_{\delta\rightarrow 0}\lim_{N\rightarrow \infty}
 \int_{{\mathbb R}^n} f_2(x)\Big[\int_{\delta}^N (Q^{\ast}_{t^m})^2({\mathcal
{I}}-P^{\ast}_{t^m})({\mathcal {I}}-P^{\ast}_{r_B^m})g(x){dt\over t}\Big]dx \nonumber\\
&=&{\rm I+II},
\end{eqnarray}
where    $f_1=f\chi_{{4B}}$, $f_2=f\chi_{(4B)^c}$ and $\chi_E$ stands for the characteristic
function of $E\subseteq{\mathbb R}^n$.

We first consider the term I. Since  $g\in {\rm L}^q(B)$, where $q=p/(p-1)$,  we conclude
$({\mathcal {I}}-P^{\ast}_{r_B^m})g\in {\rm L}^q$. By Lemma \ref{lemma
4.4} (ii), we obtain
$$
({\mathcal {I}}-P^{\ast}_{r_B^m})g=\lim_{\delta\rightarrow 0}
\lim_{N\rightarrow \infty} b_m\int_{\delta}^N
(Q^{\ast}_{t^m})^2 ({\mathcal {I}}-P^{\ast}_{t^m})
({\mathcal {I}}-P^{\ast}_{r_B^m})(g){dt\over t}
$$
 in ${\rm L}^q$.
Hence
\begin{eqnarray*}
{\rm I}&=&\lim_{\delta\rightarrow 0}\lim_{N\rightarrow \infty}
\int_{{\mathbb R}^n} f_1(x)
 \Big[ \int_{\delta}^N
(Q^{\ast}_{t^m})^2 ({\mathcal {I}}-P^{\ast}_{t^m})
({\mathcal {I}}-P^{\ast}_{r_B^m})(g)(x){dt\over t}  \Big] dx \\
 &=&  b_m^{-1}\int_{{\mathbb R}^n}  f_1(x)  ({\mathcal {I}}-P^{\ast}_{r_B^m})g(x)dx.
\end{eqnarray*}

In order to estimate  the term II, we need to   show that
 for all $y\not\in 4B$, there exists a constant $c=c(g, L)$ such that
\begin{eqnarray}
 \sup_{\delta>0,\ \!N>0}\Big|\int_{\delta}^N
(Q^{\ast}_{t^m})^2 ({\mathcal {I}}-P^{\ast}_{t^m})
({\mathcal {I}}-P^{\ast}_{r_B^m})g(x) {dt\over t}\Big|
\leq c(1+|x-x_0|)^{-(n+\epsilon)}.
\label{e4.10}
\end{eqnarray}
 To this end, set
$$
\Psi_{t,s} (L^{\ast})h(y)=(2t^m+s^m)^3\Big({d^3P^{\ast}_r\over dr^3}
\Big|_{r=2t^m+s^m} ({\mathcal I}-P^{\ast}_{t^m})h\Big)(y).
$$
Note that
$$
({\mathcal {I}}-P^{\ast}_{r_B^m})g=m\int_0^{r_B}Q^{\ast}_{s^m}(g){s^{-1}ds}.
$$
So, we use (\ref{e2.3}) and (\ref{e2.4}) to deduce
\begin{eqnarray*}
\Big|\int_{\delta}^N &&\hspace{-1cm}
(Q^{\ast}_{t^m})^2 ({\mathcal {I}}-P^{\ast}_{t^m})
({\mathcal {I}}-P^{\ast}_{r_B^m})g(x){dt\over t}\Big|\\
&=&\Big|\int_{\delta}^N  \int_0^{r_B}
(Q^{\ast}_{t^m})^2 Q^{\ast}_{s^m} ({\mathcal {I}}-P^{\ast}_{t^m})g(x){ds
\over s}{dt \over t} \Big|\\
 &\leq& c\int_{\delta}^N \int_0^{r_B}
{t^{2m}s^m\over (t^m+s^m)^3} |\Psi_{t,s} (L)g(x)|{ds
\over s}{dt \over t} \\
&\leq&c \int_{\delta}^N\int_0^{r_B} \int_{B(x_0,r_B)}
{t^{2m}s^m\over (t^m+s^m)^3}
{(t+s)^{\epsilon}\over (t+s+|x-y|)^{n+\epsilon}}|g(y)|
{dyds\over s} {dt \over t}.
\end{eqnarray*}
Because $x\not\in 4B$ yields $|x-y|\geq |x-x_0|/2$, the inequality
$$
{t^{2m}s^m(t+s)^{\epsilon}\over (t^m+s^m)^3}
\leq c\ \! {\rm min} \Big(
(ts)^{\epsilon/2}, t^{-\epsilon/2}s^{3\epsilon/2}\Big),
$$
together with H\"older's inequality and elementary integration, produces a positive constant
$c$ independent of $\delta, N>0$ such that  for all $x\not\in 4B$,
\begin{eqnarray*}
\Big|\int_{\delta}^N
Q^2_{t^m} ({\mathcal {I}}-P_{t^m})g(y){dt\over t}\Big|
 &\leq&  cr_B^{\epsilon} |x-x_0|^{-(n+\epsilon)}\|g\|_{{\rm L}^1}\\
 &\leq&  cr_B^{\epsilon+{n\over 2}}\|g\|_{{\rm L}^2}
|x-x_0|^{-(n+\epsilon)}
\end{eqnarray*}
Accordingly, (\ref{e4.10}) follows readily.

We now estimate the term II. For $f \in
 {\mathcal M}_{p}$, we derive $f\in {\rm L}^p\big((1+|x|)^{-(n+\epsilon_0)}dx\big)
 $.
 The estimate
(\ref{e4.10}) yields a constant $c>0$ such that
\begin{eqnarray*}
\sup_{\delta>0,\ \!N>0}\int_{{\mathbb R}^n} \Big|f_2(x)
 \int_{\delta}^N
(Q^{\ast}_{t^m})^2 ({\mathcal {I}}-P^{\ast}_{t^m})
({\mathcal {I}}-P^{\ast}_{r_B^m})(g)(x){dt\over t}\Big|  dx
\leq c.
\end{eqnarray*}
This allows us to
pass  the limit inside the integral  of  II. Hence
\begin{eqnarray*}
{\rm II}&=&\lim_{\delta\rightarrow 0}\lim_{N\rightarrow \infty}
\int_{{\mathbb R}^n} f_2(x)
  \Big[\int_{\delta}^N
(Q^{\ast}_{t^m})^2 ({\mathcal {I}}-P^{\ast}_{t^m})
({\mathcal {I}}-P^{\ast}_{r_B^m})(g)(x){dt\over t} \Big]  dx \\
&=&\int_{{\mathbb R}^n} f_2(x)\Big(
\lim_{\delta\rightarrow 0}\lim_{N\rightarrow \infty}
\Big[\int_{\delta}^N
(Q^{\ast}_{t^m})^2 ({\mathcal {I}}-P^{\ast}_{t^m})({\mathcal {I}}-P^{\ast}_{r_B^m})
(g)(x){dt\over t} \Big]\Big)  dx \\
 &=&  b_m^{-1}\int_{{\mathbb R}^n}  f_2(x)  ({\mathcal {I}}-P^{\ast}_{r_B^m})g(x)dx.
\end{eqnarray*}
Combining the previous formulas for I and II, we obtain the identity (\ref{e4.8}).
\end{proof}

\begin{remark} For a background of Proposition \ref{prop 4.3}, see also \cite[Proposition 5.1]{DY2}.
\end{remark}

\section{Description through Littlewood-Paley function}
\setcounter{equation}{0}

\subsection{The space \LV\ as the dual of the atomic space} Following \cite{Z}, we give the
 following definition.

\begin{definition} Let $1<p<\infty$, $q=p/(p-1)$ and $\lambda\in (0,n)$. Then

\item{\rm(i)} A complex-valued function $a$ on $\mathbb R^n$ is called a $(q,\lambda)$-atom provided:

\ {($\alpha$)} $a$ is supported on an open ball $B\subset{\mathbb R}^n$ with radius $r_B$;

\ {($\beta$)} $\int_{\mathbb R^n} a(x)dx=0$;

\ {($\gamma$)} $\|a\|_{{\rm L}^q}\leq r_B^{-\lambda/p}$.

\item{\rm(ii)} ${\rm H}^{q,\lambda}({\mathbb R}^n)$ comprises those linear
functionals admitting an atomic decomposition
$f=\sum_{j=1}^{\infty}\eta_ja_j,$ where $a_j$'s are  $(q, \lambda)$-atoms, and
$\sum_{j}|\eta_j|<\infty.$
\end{definition}

The forthcoming result reveals that ${\rm H}^{q,\lambda}({\mathbb R}^n)$ exists as
a predual of ${\rm L}^{p,
\lambda}({\mathbb R}^n)$.

\begin{prop}\label{prop4.1}
Let $1<p<\infty$, $q=p/(p-1)$ and $\lambda\in (0,n)$. Then ${\rm L}^{p,
\lambda}({\mathbb R}^n)$ is the dual $\big({\rm H}^{q,\lambda}({\mathbb
R}^n)\big)^\ast$ of ${\rm H}^{q,\lambda}({\mathbb
R}^n)$. More precisely, if
$h=\sum_j\eta_ja_j\in {\rm H}^{q,\lambda}({\mathbb R}^n)$ then
$$
\langle h, \ell\rangle=\lim_{k\rightarrow\infty}\sum_{j=1}^{k}\eta_j\int_{{\mathbb R}^n}
a_j(x)\ell(x)dx
$$
is a well-defined continuous linear functional for each $\ell\in {\rm L}^{p, \lambda}({\mathbb
R}^n)$, whose norm is equivalent to $\|\ell\|_{{\rm L}^{p,
\lambda}}$; moreover, each continuous linear functional on
${\rm H}^{q,\lambda}({\mathbb R}^n)$ has this form.
\end{prop}

\begin{proof}
See \cite[Proposition 5]{Z} for a proof of Proposition \ref{prop4.1}.
\end{proof}

\subsection{Characterization of \LV\ by means of Littlewood-Paley function} We now state
a full characterization of ${\rm L}^{p, \lambda}({\mathbb R}^n)$
 space for  $1<p<\infty$ and $\lambda\in (0, n).$ For the case $p=2$, see also \cite[Lemma 2.1]{X}
  as well as \cite[Theorem 1 (i)] {WX}.

\begin{prop} \label{prop4.5a} Let $1<p<\infty$,
$\lambda\in (0, n)$ and $f\in {\mathcal M}_{\sqrt{\triangle}, p}$. Then the following two conditions are
equivalent:

\item{\rm (i)} $f\in {\rm L}^{p, \lambda}({\mathbb
R}^n)$;

\item{\rm (ii)}
\begin{eqnarray*}
I(f,p)=\sup_{B\subset{\mathbb
R}^n}r_B^{-\frac{\lambda}{p}}\Big{\|}\Big\{\int_0^{r_B}|t
{\partial\over
\partial t} e^{-t\sqrt{\triangle}}f(x)|^2{dt\over
t}\Big\}^{1/2}\Big{\|}_{{\rm L}^p(B)}<\infty,
\end{eqnarray*}
where the supremum is taken over all Euclidean open ball
$B\subset{\mathbb R}^n$ with radius $r_B$.
\end{prop}

\begin{proof} It suffices to verify (ii)$\Rightarrow$(i) for which the reverse implication
follows readily from \cite[Theorem 2.1]{FJN}. Suppose (ii) holds.
Proposition \ref{prop4.1} suggests us to show $f\in
\big({\rm H}^{\frac{p}{p-1},\lambda}({\mathbb R}^n\big)^\ast$ in order to verify (i).
Now, let $g$ be a $(\frac{p}{p-1},\lambda)$-atom and
$$
p_t(x)=\frac{c_nt}{(t^2+|x|^2)^{\frac{n+1}{2}}}.
$$
Then for any open ball $B\subset{\mathbb R}^n$ with radius $r_B$ and
its tent
$$
T(B)=\{(x,t)\in{\mathbb R}^{n+1}_+:\quad x\in B,\ t\in (0,r_B)\},
$$
we have (cf. \cite[p.183]{St})
\begin{eqnarray*}
|\langle f,g\rangle|&=&\Big|\int_{{\mathbb R}^n}f(x)g(x)dx\Big|\\
&=&4\left|\int_{{\mathbb R}^n}\int_0^\infty\Big(t\frac{\partial}{\partial t}p_t\ast
f(x)\Big)\Big(t\frac{\partial}{\partial t}p_t\ast g(x)\Big)t^{-1}dtdx\right|\\
&\leq& 4\big(I(B)+J(B)\big).
\end{eqnarray*}
Here,
\begin{eqnarray*}
I(B)&=&\int_{4B}\int_0^{r_{4B}}\Big|t\frac{\partial}{\partial t}p_t\ast
f(x)\Big|\Big|t\frac{\partial}{\partial t}p_t\ast g(x)\Big|t^{-1}dtdx\\
\\
&\leq&\left(\int_{4B}\Big(\int_0^{r_{4B}}\Big|t\frac{\partial}{\partial
t}p_t\ast
f(x)\Big|^2t^{-1}dt\Big)^\frac{p}{2}dx\right)^\frac1p\\
&&\hspace{1.6cm}\times\left(\int_{4B}\Big(\int_0^{r_{4B}}\Big|t\frac{\partial}{\partial
t}p_t\ast
g(x)\Big|^2t^{-1}dt\Big)^\frac{p}{2(p-1)}dx\right)^\frac{p-1}{p}\\
&\leq& c r_B^\frac{\lambda}{p}I(f,p)\|g\|_{{\rm L}^\frac{p}{p-1}({\mathbb
R}^n)}\\
&\leq& c I(f,p),
\end{eqnarray*}
due to H\"older's inequality, the ${\rm L}^\frac{p}{p-1}$-boundedness of
the Littlewood-Paley ${\mathcal G}$-function, and $g$ being a
$(\frac{p}{p-1},\lambda)$-atom.

Meanwhile,
\begin{eqnarray*}
J(B)&=&\sum_{k=1}^\infty\int_{T(2^{k+1}B)\setminus T(2^kB)}
\Big|t\frac{\partial}{\partial t}p_t\ast f(x)\Big|\Big|t\frac{\partial}{\partial t}p_t\ast g(x)
\Big|t^{-1}dtdx\\
\\
&\leq&c\sum_{k=1}^\infty\Big\|\Big\{\int_0^{2^{k+1}r_B}\Big|t\frac{\partial}{\partial
t}p_t\ast
f(x)\Big|^2t^{-1}dt\Big\}^\frac12\Big\|_{{\rm L}^p(2^{k+1}B)}\\
&&\times\Big\|\Big\{\int_0^{2^{k+1}r_B}\Big|t\frac{\partial}{\partial
t}p_t\ast g(x)\Big|^2t^{-1}dt\Big\}^\frac12\Big\|_{{\rm L}^\frac{p}{p-1}(2^{k+1}B)}\\
&\leq&c\sum_{k=1}^\infty(2^kr_B)^\frac{\lambda}{p}I(f,p)2^{-\frac{kn}{p}}r_B^{-\frac{\lambda}{p}}\\
&\leq&c I(f,p),
\end{eqnarray*}
for which we have used the H\"older inequality and the fact that if
$|y-x|\ge 2^k r_B$ then
$$
\Big|t\frac{\partial}{\partial t}p_t\ast g(x)\Big|\leq c
\frac{t^3}{(2^kr_B)^{3+n}}\|g\|_{{\rm L}^1(B)}\leq
c\frac{t^3}{(2^kr_B)^{3+n}}r_B^\frac{n-\lambda}{p}
$$
for the $(\frac{p}{p-1},\lambda)$-atom $g$. Accordingly, $f\in$\ \LV.
\end{proof}

\subsection{Characterization of \LNV\ by means of Littlewood-Paley function} Of course,
it is natural to explore a characterization of ${\rm
L}_{L}^{p,\lambda}({\mathbb R}^n)$ similar to Proposition \ref{prop4.5a}.

\begin{prop} \label{prop 4.2} Let $1<p<\infty$, $\lambda\in (0,n)$ and $f\in {\mathcal M}_{p}$.
Assume that  $L$ satisfies the assumptions   (a) and  (b) of Section
2.2. Then the following two conditions are equivalent:

\item{\rm(i)} $f\in$\ \LNV ;

\item{\rm (ii)}

\begin{eqnarray*}
|\!|\!|f|\!|\!|_{{\rm L}_L^{p,\lambda}}=\sup_{B\subset{\mathbb
R}^n}r_B^{-\frac\lambda{p}}\Big{\|}\Big\{\int_0^{r_B}|Q_{t^m}
({\mathcal {I}}-P_{t^m})f(x)|^2{dt\over
t}\Big\}^{1/2}\Big{\|}_{{\rm L}^p(B)}<\infty,
\end{eqnarray*}
where the supremum is taken over all Euclidean open ball
$B\subset{\mathbb R}^n$ with radius $r_B$.
\end{prop}

\begin{proof} (i)$\Rightarrow$(ii). Suppose $f\in$\ \LNV.\ Note that
$$
Q_{t^m}({\mathcal {I}}-P_{t^m})=Q_{t^m}({\mathcal {I}} -P_{t^m})
({\mathcal {I}}-P_{r_{B}^m})+ Q_{t^m}(I-P_{t^m})P_{r_{B}^m}.
$$
So, we turn to verify both
\begin{equation}
\Big{\|}\Big\{\int_0^{r_B}|Q_{t^m}({\mathcal {I}} -P_{t^m})
({\mathcal {I}}-P_{r_{B}^m})f(x)|^2{dt\over t}\Big\}^{1/2}\Big{\|}_{{\rm L}^p(B)}
\leq cr_B^{\lambda \over p}\|f\|_{{\rm L}_L^{p, \lambda}}
\label{e4.3}
\end{equation}
and
\begin{equation}
\Big{\|}\Big\{\int_0^{r_B}|Q_{t^m}(I-P_{t^m})P_{r_{B}^m}
f(x)|^2{dt\over t}\Big\}^{1/2}\Big{\|}_{{\rm L}^p(B)}
\leq cr_B^{\lambda \over p}\|f\|_{{\rm L}_L^{p, \lambda}},
\label{e4.4}
\end{equation}
thereby proving (ii). To do so, we will adapt the argument on
 pp. 85-86 of \cite{J} to present situation -- see also page 955 of \cite{DY2}. To prove
(\ref{e4.3}),  let us consider  the square function
${\mathcal G}(h)$  given  by
$$
{\mathcal G}(h)(x)=\Big(\int_0^{\infty}
|Q_{t^m}({\mathcal {I}}-P_{t^m})h(x)|^2{dt\over t}\Big)^{1/2}.
$$
From (\ref{e2.7}), the function ${\mathcal G}(h)$ is bounded on
${\rm L}^p({\mathbb R}^n)$ for   $1<p<\infty$. Let $b=b_1+b_2$, where
$b_1=({\mathcal I}-P_{r^m_{B}})f\chi_{2B}$, and $b_2=({\mathcal
I}-P_{r^m_{B} })f\chi_{(2B)^c}$. Using Lemma \ref{e3.3}, we obtain
\begin{eqnarray}
\Big{\|}&&\hspace{-1cm}\Big\{\int_0^{r_B}|Q_{t^m}({\mathcal {I}} -P_{t^m})
b_1(x)|^2{dt\over t}\Big\}^{1/2}\Big{\|}_{{\rm L}^p(B)}\nonumber\\
&\leq&\Big{\|}\Big\{\int_0^{\infty}|Q_{t^m}({\mathcal {I}} -P_{t^m})
b_1(x)|^2{dt\over t}
\Big\}^{1/2}\Big{\|}_{{\rm L}^p}\nonumber\\
&\leq&c\|{\mathcal G}(b_1)\|_{{\rm L}^p}\nonumber\\
&\leq& c\|b_1\|_{{\rm L}^p} \nonumber\\
&=&c\Big(\int_{2B}|({\mathcal I}-P_{r^m_{B} })f(x)|^pdx\Big)^{1/p}\nonumber\\
&\leq& c\Big(\int_{2B}|({\mathcal I}-P_{ r^m_{2B} })f(x)|^pdx\Big)^{1/p} +
cr_B^{n/p}\cdot \sup_{x\in{2B}}|P_{r^m_{B}}f(x)-P_{r^m_{2B
} }f(x)|^p\nonumber\\
&\leq& cr_B^{\lambda \over p}\|f\|_{{\rm L}_L^{p, \lambda}}.
\label{e4.5}
\end{eqnarray}
On the other hand, for any $x\in {B}$ and $y\in (2B)^c$, one has
$|x-y|\geq r_B$. From  Proposition 3.5, we obtain
\begin{eqnarray*}
|Q_{t^m }({\mathcal {I}}-P_{ t^m })b_2(x)|
&\leq&c\int_{{\mathbb R}^n\backslash  2B}{{t^{\epsilon}}\over (t+
|x-y|)^{n+{\epsilon}}}
|({\mathcal I}-P_{r^m_{{B}} })f(y)|dy\\
&\leq&c\Big({t\over r_{B}}\Big)^{\epsilon}\int_{{\mathbb R}^n}
{{r^{\epsilon}_{B}}\over
(r_{B}+|x-y|)^{n+{\epsilon}}} |({\mathcal I}-P_{r^m_{B}})f(y)|dy\\
&\leq&c\Big({t\over r_{B}}\Big)^{\epsilon}
r_B^{\lambda-n\over p}\|f\|_{{\rm L}_L^{p, \lambda}},
\end{eqnarray*}
which implies
\begin{eqnarray*}
\Big{\|} \Big\{\int_0^{r_B}|Q_{t^m}({\mathcal {I}} -P_{t^m})
b_2(x)|^2{dt\over t}\Big\}^{1/2}\Big{\|}_{{\rm L}^p(B)}
 \leq  cr_B^{\lambda \over p}\|f\|_{{\rm L}_L^{p, \lambda}}.
\end{eqnarray*}
This, together with    (\ref{e4.5}), gives (\ref{e4.3}).

Next, let us check (\ref{e4.4}). This time, we have $0<t< r_B$, whence getting
from Lemma \ref{e3.3} that for any $x\in{\mathbb R}^n,$
$$
|P_{{1\over 2}r^m_{{B}} }f(x)-P_{{(t^m+{1\over 2}r^m_{{B}})}}f(x)|
\leq cr_B^{\lambda-n\over p}\|f\|_{{\rm L}_L^{p, \lambda}}.
$$
By (\ref{e2.4}), the kernel $K_{t, r_{B}}(x,y)$ of the
operator
$$
Q_{t^m}P_{{1\over 2}r^m_{B}}={t^m\over t^m+{1\over 2}r^m_{B}}
Q_{{(t^m+{1\over 2}r^m_{B})}}
$$
satisfies
$$
|K_{t, r_{B}}(x,y)|\leq c \Big({t\over r_{B}}\Big)^m{{r_{B}^{\epsilon}}\over
(r_{B} +|x-y|)^{n+{\epsilon}}}.
$$
Using the commutative property of the semigroup $\{e^{-tL}\}_{t>0}$ and the
estimate (\ref{e2.4}), we deduce
\begin{eqnarray*}
|Q_{t^m}(I-P_{t^m})P_{r_{B}^m}
f(x)|
&=&|Q_{t^m}P_{{1\over 2}r^m_{B}} (P_{{1\over 2}r^m_{B} }-
P_{{(t^m+{1\over 2}r^m_{B})}})f(x)|\\
&\leq&c\Big({t\over r_{B}}\Big)^m\int_{{\mathbb R}^n}{{r_{B}^{\epsilon}}\over
(r_{B} +|x-y|)^{n+{\epsilon}}}|
(P_{{1\over 2}r^m_{B} }-P_{{(t^m+{1\over 2}r^m_{B})}})f(y)|dy\\
&\leq&c\Big({t\over r_B}\Big)^mr_B^{\lambda-n \over p}
\|f\|_{{\rm L}_L^{p, \lambda}},
\end{eqnarray*}
whence deriving
\begin{eqnarray*}
\Big{\|}\Big\{\int_0^{r_B}|Q_{t^m}(I-P_{t^m})P_{r_{B}^m}
f(x)|^2{dt\over t}\Big\}^{1/2}\Big{\|}_{{\rm L}^p(B)}\leq
cr_B^{\lambda \over p}\|f\|_{{\rm L}_L^{p, \lambda}}.
\end{eqnarray*}
This gives (\ref{e4.4}) and consequently (ii).

(ii)$\Rightarrow$ (i). Suppose (ii) holds. The duality argument
for ${\rm L}^p$ shows that for any open ball $B\subset{\mathbb R}^n$ with radius $r_B$,
\begin{eqnarray}
\Big(r_B^{-\lambda}\int_B |f(x)-P_{r_B^m}f(x)|^pdx\Big)^{1/p}
&=&\sup\limits_{\|g\|_{{\rm L}^q(B)\leq 1}}r_B^{-\lambda/p}\Big|\int_{{\mathbb R}^n}
(I-P_{r_B^m})f(x)g(x)dx\Big|\nonumber\\
&=&\sup\limits_{\|g\|_{{\rm L}^q(B)\leq 1}}r_B^{-\lambda/p}\Big|\int_{{\mathbb R}^n}
f(x)(I-P^{\ast}_{r_B^m})g(x)dx\Big|.
\label{e4.11}
\end{eqnarray}
Using the identity (\ref{e4.8}), the estimate
(\ref{e4.7}) and the H\"older inequality, we have
\begin{eqnarray}
 \Big|\int_{{\mathbb R}^n}
f(x)(I-P^{\ast}_{r_B^m})g(x)dx\Big|
 &\leq& c\int_{{\mathbb R}^{n+1}_+}|Q_{t^m} ({\mathcal {I}}-P_{t^m})f(x)  \
 Q^{\ast}_{t^m} ({\mathcal {I}}-P^{\ast}_{r_B^m})g(x)|{dxdt\over t}\nonumber \\
 &\leq&cr_B^{\lambda/p}|\!|\!|f|\!|\!|_{{\rm L}_L^{p,\lambda}}\|g\|_{{\rm L}^q}.
\label{e4.12}
\end{eqnarray}
Substituting (\ref{e4.12}) back to (\ref{e4.11}), by Definition 3.1
we find a constant $c>0$ such that
$$
\|f\|_{{\rm L}_L^{p, \lambda}}\leq c|\!|\!|f|\!|\!|_{{\rm L}_L^{p,\lambda}}<\infty.
$$
This just proves $f\in {\rm L}_L^{p,
\lambda}({\mathbb R}^n)$, thereby yielding (i).
\end{proof}

\begin{remark} In the case of $p=2$, we can interpret Proposition \ref{prop 4.2} as a
measure-theoretic characterization, namely, $f\in {\rm L}_L^{2,\lambda}({\mathbb R}^n)$
when and only when
$$
d\mu_f(x,t)=|Q_{t^m}
({\mathcal {I}}-P_{t^m})f(x)|^2{dxdt\over
t}
$$
is a $\lambda$-Carleson measure on ${\mathbb R}^{n+1}_+$. According to
\cite[Lemma 4.1]{EJPX}, we find further that $f\in {\rm L}_L^{2,\lambda}({\mathbb R}^n)$
is equivalent to
$$
\sup_{(y,s)\in{\mathbb R}^{n+1}_+}\int_{{\mathbb R}^{n+1}_+}
\left(\frac{s}{\big(|x-y|^2+(t+s)^2\big)^\frac{n+1}{2}}\right)^\lambda d\mu_f(x,t)<\infty.
$$
\end{remark}

\subsection{A sufficient condition for \LNV\ = \LV}
In what follows, we assume that $L$ is a linear operator of type $\omega$ on
${\rm L}^2({\mathbb R}^n)$ with $\omega<\pi/2$ -- hence $L$ generates
  an analytic semigroup
$e^{-zL}, 0\leq |{\rm Arg}(z)|<\pi/2-\omega$. We also assume that for
each
 $t>0$, the kernel  $p_t(x,y)$
of  $ e^{-tL} $  is H\"older continuous in both variables
$x$, $y$  and
   there exist positive
constants $m $, $\beta>0 $ and $0<\gamma\leq 1
$ such that for all $t>0$, and
$x,y, h\in{\mathbb R}^n$,
\begin{eqnarray}
|p_t(x,y)|\leq c {t^{\beta/m}\over  (t^{1/m}+|x-y|)^{n+\beta}}\quad\forall\ t>0,\ x,y\in{\mathbb R}^n,
\label{e5.1}
\end{eqnarray}
\begin{eqnarray}
\label{e5.2}
&{}&|p_t(x+h,y)-p_t(x,y)| + |p_t(x,y+h)-p_t(x,y)|\nonumber\\
 &\leq& c |h|^{\gamma}{t^{\beta/m}\over  (t^{1/m}+|x-y|)^{n+\beta+\gamma}}\quad\forall h\in{\mathbb R}^n
\quad\hbox{with}\quad 2|h|\leq t^{1/m}+|x-y|,
\end{eqnarray}
and
\begin{eqnarray}
\int_{{\mathbb R}^n}p_t(x,y)dx=
\int_{{\mathbb R}^n}p_t(x,y)dy=1\quad\forall t>0.
\label{e5.3}
\end{eqnarray}

\medskip

\begin{prop} \label{theorem 5.2} Let $1<p<\infty$ and $\lambda\in (0,n)$. Given
 an operator $L$ which generates a semigroup $e^{-tL}$ with the heat kernel bounds
 (\ref{e2.2}) and (\ref{e2.3}). Assume that  $L$ satisfies the conditions
 (\ref{e5.1}), (\ref{e5.2}) and (\ref{e5.3}).
Then \LNV\ and\ \LV\ coincide, and their
norms are equivalent.
\end{prop}

\begin{proof} Since Proposition \ref{prop3.2} tells us that \LV\ $\subseteq$ \LNV\
under the above-given conditions, we only need to check \LNV $\subseteq$ \LV.
Note that \LV \ is the dual of ${\rm H}^{q,\lambda}({\mathbb R}^n)$, $q=p/(p-1)$.  It
reduces to prove that if $f\in$ \LNV, then $f\in
({\rm H}^{q,\lambda}({\mathbb R}^n))^{\ast}$.  Let $g$ be a $(q, \lambda)$-atom. Using the conditions
(\ref{e5.1}), (\ref{e5.2}) and  (\ref{e5.3}) of the operator $L$,
together with the properties of
  of $(q, \lambda)$-atom of $g$,  we can follow the argument for Lemma \ref{lemma 4.4} (ii) to verify
\begin{eqnarray*}
\int_{{\mathbb R}^n}f(x)g(x)dx
=b_m\int_{{\mathbb R}^{n+1}_+}Q_{t^m} ({\mathcal {I}}-P_{t^m})f(x)
Q_{t^m}^{\ast} g(x){dxdt\over t}\quad\hbox{where}\quad b_m=\frac{36m}{5}.
\end{eqnarray*}
Consequently,
\begin{eqnarray*}
|\langle f,g\rangle|&=&\Big|\int_{{\mathbb R}^n}f(x)g(x)dx\Big|\\
&=&\Big|\int_{{\mathbb R}^{n+1}_+}Q_{t^m} ({\mathcal {I}}-P_{t^m})f(x)
Q_{t^m}^{\ast} g(x){dxdt\over t}\Big|\\
&\leq&\int_{T(4B)} \big|Q_{t^m} ({\mathcal {I}}-P_{t^m})f(x)
Q_{t^m}^{\ast} g(x)\big|{dxdt\over t}\\
&& +\sum_{k=1}^{\infty}\int_{T(2^{k+1}B)\backslash T(2^{k}B) }
\big|Q_{t^m} ({\mathcal {I}}-P_{t^m})f(x)\ Q_{t^m}^{\ast}g(x)\big|{dxdt\over t}\\
&=& {\rm D_1} +\sum_{k=2}^{\infty} {\rm D_k}.
\end{eqnarray*}
Define  the Littlewood-Paley function ${\mathcal G} h$  by
$$
{\mathcal G}(h)(x)=\bigg[ \int_0^{\infty} |Q_{t^m}^{\ast} h(x)|^2{dt\over t}\bigg]^{1/2}.
$$
By (\ref{e2.7}), ${\mathcal G}$  is   bounded on
${\rm L}^p({\mathbb R}^n)$ for  $1<p<\infty$.

Following the proof of Lemma \ref{lemma 4.4} (i), together with the
property ($\gamma$) of $(q, \lambda)$-atom $g$, we derive
\begin{eqnarray*}
{\rm D_1}
&\leq&\Big\|\Big\{\int_0^{r_{2B}}\big|Q_{t^m} ({\mathcal {I}}-P_{t^m})f(x)
\big|^2{dt\over t}\Big{\}}^{1/2}\Big\|_{{\rm L}^p(2B)}
\Big\|\Big\{\int_0^{r_{2B}}\big|Q_{t^m}^{\ast}g(x)\big|^2{dt\over t}\Big{\}}^{1/2}\Big\|_{{\rm L}^q(2B)}
\\
&\leq&\Big\|\Big\{\int_0^{r_{2B}}\big|Q_{t^m} ({\mathcal {I}}-P_{t^m})
f(x)\big|^2{dt\over t}\Big{\}}^{1/2}\Big\|_{{\rm L}^p(2B)}
\|{\mathcal G}(g)\|_{{\rm L}^q}\\
&\leq&cr_B^{\lambda \over p}|\!|\!|f|\!|\!|_{{\rm L}_{L}^{p, \lambda}}\|g\|_{{\rm L}^q}
 \leq c\|f\|_{{\rm L}_{L}^{p, \lambda}}.
\end{eqnarray*}
On the other hand, we note that for   $x\in T(2^{k+1}B)\backslash T(2^{k}B)$
and $y\in B$, we have that $|x-y|\geq 2^kr_B$. Using the estimate (\ref{e2.4}) and
the properties ($(\alpha$) and ($\gamma$) of $(q, \lambda)$-atom $g$,
 we obtain
 \begin{eqnarray*}
|Q_{t^m}^{\ast}g(x)|&\leq&c\int_{B} {t^{\epsilon}
\over (t+ |x-y|)^{n+{\epsilon} }}|g(y)|dy\\
 &\leq& c  {t^{\epsilon}
\over (2^kr_B)^{n+{\epsilon} }}\int_{B}|g(y)|dy\\
&\leq&c  {t^{\epsilon}
\over (2^kr_B)^{n+{\epsilon} }}r_B^{n-\lambda\over p},
\end{eqnarray*}
which implies
\begin{eqnarray*}
\Big{\|}\Big\{\int_0^{2^kr_B}
|Q_{t^m}^{\ast}g(x)\chi_{T(2^{k+1}B)\backslash T(2^{k}B) }|^2
{dt\over t}\Big\}^{1/2}\Big{\|}_{{\rm L}^{q}(2^kB)}
\leq c2^{kn( {1\over q}-1)}r_B^{-{\lambda\over p}}
\end{eqnarray*}
Therefore,
\begin{eqnarray*}
{\rm D}_k
&\leq&\Big{\|}\Big\{\int_0^{2^kr_B}|Q_{t^m} ({\mathcal {I}}-P_{t^m}) f(x)
|^2{dt\over t}\Big\}^{1/2}\Big{\|}_{{\rm L}^p(2^kB)}\\
&&\times\Big{\|}\Big\{\int_0^{2^kr_B}|Q_{t^m}^{\ast}g(x)\chi_{T(2^{k+1}B)\backslash
T(2^{k}B) }|^2
{dt\over t}\Big\}^{1/2}\Big{\|}_{{\rm L}^{q}(2^kB)}\\
&\leq&c (2^kr_B)^{\lambda \over p}2^{kn({1\over q}-1)}r_B^{-{\lambda  \over p}}
|\!|\!|f|\!|\!|_{{\rm L}_{L}^{p, \lambda}}\\
&\leq&c2^{k(\lambda- n)\over p}
\|f\|_{{\rm L}_{L}^{p, \lambda}}.
\end{eqnarray*}
Since $\lambda\in (0, n)$, we have
\begin{eqnarray*}
|\langle f, g\rangle|
\leq c\|f\|_{{\rm L}_{L}^{p, \lambda}}+
c\sum_{k=1}^{\infty} 2^{k(\lambda-n)\over p}
\|f\|_{{\rm L}_{L}^{p, \lambda}}
\leq c\|f\|_{{\rm L}_{L}^{p, \lambda}}.
\end{eqnarray*}
This, together with Proposition \ref{prop4.1}, implies
  $f\in \big({\rm H}^{q,\lambda}({\mathbb R}^n)\big)^{\ast}=$ \LV\ .
\end{proof}

\bigskip

{\footnotesize { DEPARTMENT OF MATHEMATICS, MACQUARIE UNIVERSITY, NSW 2109, AUSTRALIA}

{\it E-mail address}: duong@ics.mq.edu.au

\bigskip

DEPARTMENT OF MATHEMATICS AND STATISTICS, MEMORIAL UNIVERSITY
OF NEWFOUNDLAND, ST. JOHN'S, NL, A1C 5S7, CANADA

{\it E-mail address}: jxiao@math.mun.ca
\bigskip

DEPARTMENT OF MATHEMATICS, ZHONGSHAN  UNIVERSITY, GUANGZHOU, 510275, P.R. CHINA

{\it E-mail address}: mcsylx@mail.sysu.edu.cn

}

\end{document}